\documentclass[a4paper,12pt]{article}
\usepackage{amsmath,amsthm,amsfonts,amssymb}
\usepackage{graphicx}
\usepackage{enumerate}
\usepackage[longnamesfirst]{natbib}
\bibpunct[ ]{(}{)}{,}{a}{}{,}
\usepackage[hmargin=1in,vmargin=1in]{geometry}
\usepackage[dvipsnames]{xcolor}
\definecolor{DarkGreen}{rgb}{0,0.35,0}
\usepackage{hyperref}
\usepackage{pgfplots}
\pgfplotsset{compat=1.18}
\hypersetup{
    colorlinks=true,
    linkcolor=NavyBlue,    
    urlcolor=TealBlue,     
    citecolor=OliveGreen,  
    filecolor=Maroon       
}

\renewcommand{\d}{\mathrm d}
\newcommand{\R}{\mathbb R}
\newcommand{\Q}{\mathbb Q}
\newcommand{\Z}{\mathbb Z}

\newcommand{\wt}{\widetilde}
\newcommand{\wh}{\widehat}
\newcommand{\ol}{\overline}

\newcommand{\eps}{\varepsilon}

\newcommand{\D}{\mathcal D}

\newcommand{\Hc}{\mathcal H}

\renewcommand{\O}{\mathcal{O}}

\newcommand{\e}{\mathrm e}

\newcommand{\RW}{\mathrm{RW}}
\newcommand{\Br}{\mathrm{Br}}
\newcommand{\SJ}{\mathrm{SJ}}
\newcommand{\LF}{\mathrm{LF}}
\newcommand{\epi}{\mathfrak e}
\newcommand{\LSC}{\mathrm{LSC}}
\newcommand{\ev}{\mathrm{ev}}

\newcommand{\cH}{\mathcal{H}}
\newcommand{\dH}{d_{\mathcal{H}}}

\newtheorem{proposition}{Proposition}[section]
\newtheorem{theorem}[proposition]{Theorem}
\newtheorem{lemma}[proposition]{Lemma}

\newtheorem{corollary}[proposition]{Corollary}
\newtheorem{conjecture}[proposition]{Conjecture}
\theoremstyle{definition}
\newtheorem*{remark}{Remark}

\newtheorem{question}[proposition]{Question}
\newtheorem{problem}[proposition]{Problem}

\numberwithin{equation}{section}

\title{The dynamical Sepp\"al\"ainen--Johansson model and the Brownian web landscape}
\author{
Julian Ransford\thanks{Department of Pure Mathematics and Mathematical Statistics, University of Cambridge, UK. E-mail: {\tt julian.ransford@gmail.com}} \and
B\'alint Vet\H o\thanks{Department of Stochastics, Institute of Mathematics,
Budapest University of Technology and Economics, M\H uegyetem rkp.\ 3., H-1111 Budapest, Hungary
and HUN--REN Alfr\'ed R\'enyi Institute of Mathematics, Re\'altanoda u.\ 13--15., H-1053 Budapest, Hungary.
E-mail: {\tt vetob@math.bme.hu}}
\and
B\'alint Vir\'ag\thanks{Departments of Mathematics and Statistics, University of Toronto, Canada. E-mail: {\tt balint@math.toronto.edu}}
}

\begin{document}

\maketitle

\begin{abstract}
The dynamical Sepp\"al\"ainen--Johansson first passage percolation defines a random directed metric on $\Z^2$
where the environment is updated according to a Poisson process.
We show that its finite dimensional distributions converge after $0:1:2$ scaling
to those of the Brownian web landscape defined on the finite dimensional marginals of the dynamical Brownian web.
We prove non-uniqueness of geodesics in the Brownian web landscape
and we relate it to exact formulas on one-sided reflected Brownian motions. We prove a locally uniform law of large numbers and uniform deviation
bounds, and identify Poisson increment-stationary laws for the
associated height evolution.
\end{abstract}

\section{Introduction}

Planar last passage percolation models are standard examples of models in the Kardar--Parisi--Zhang (KPZ) universality class.
In several integrable models, after $1:2:3$ scaling, the last passage time converges to the directed landscape, the universal scaling limit in the KPZ class.

In the present paper, we consider the Sepp\"al\"ainen--Johansson first passage percolation model first introduced in \cite{S98}
and later studied in \cite{J01}. The model is the weighted di-graph distance $F(i,j;k;\ell)$ of the points $(i,j)$ and $(k,\ell)$ in the square grid directed up/right, with edge weights zero along vertical edges and i.i.d.\ Bernoulli-$1/2$ random variables along horizontal edges.

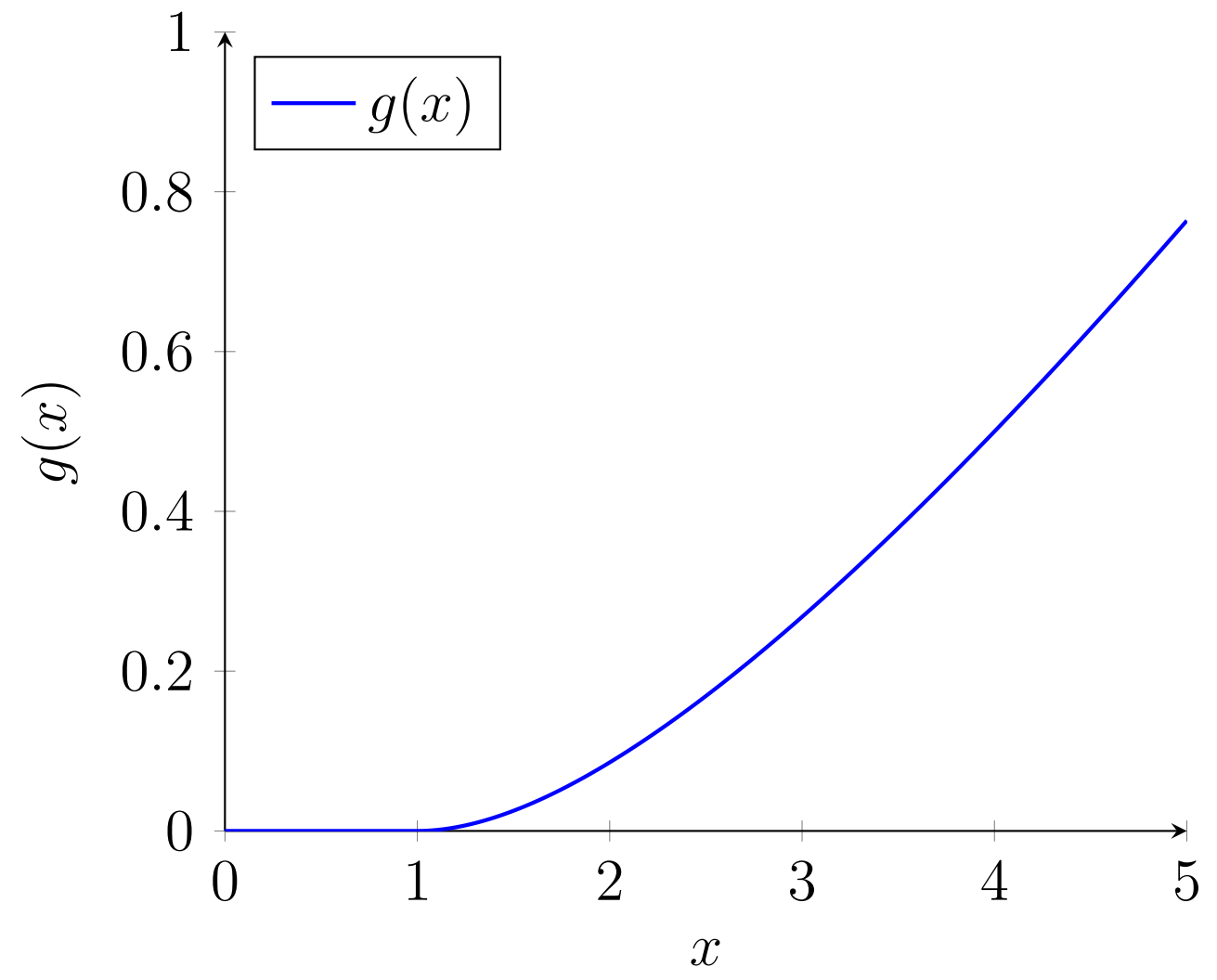
\begin{figure}
    \centering
    \begin{tikzpicture}
\begin{axis}[
  xlabel=$x$,
  ylabel=$g(x)$,
  xmin=0, xmax=5,
  ymin=0, ymax=1,
  axis lines=left,
  samples=200,
  domain=0:5,
  legend pos=north west,
]
\addplot[blue, thick] {x < 1 ? 0 : (sqrt(x) - 1)^2 / 2};
\legend{$g(x)$}
\end{axis}
\end{tikzpicture}
    \caption{The limit of the distance function in direction $(x,1)$.}
    \label{fig:limitshape}
\end{figure}

The model can also be seen as last passage percolation since the weights are bounded.
The limit shape of balls was described in \cite{S98}.
An equivalent description is given by the curve
\begin{equation}
g(x)=\lim_{\ell \to \infty}\frac{ F(0,0;\lfloor x\ell \rfloor,\ell)}{\ell}=\mathbf{1}_{x\ge 1}\frac{(\sqrt{x}-1)^2}{2}, \qquad x\ge 0,
\end{equation}
which captures the large-scale behavior of the distance function in the direction $(x,1)$; see Figure~\ref{fig:limitshape}.

In the directions $(x,1)$ with $x\in[0,1)$, we have $g(x)=0$, so the rescaled first passage values converge to zero.
For directions $(x,1)$ with $x>1$, $g$ has a positive second derivative.  In this regime, after $1:2:3$ scaling, the passage value has GUE Tracy--Widom fluctuations and the distances converge to the directed landscape, see \cite{J01}, \cite{DNV23} and \cite{DV21}.
The limit shape of a generalized Sepp\"al\"ainen--Johansson model was recently described in \cite{R24}.

The direction $(1,1)$ is critical where the curved part of the limit shape meets its flat part. This is one of the very few phenomena that have been proven to happen even in non-integrable (and non-oriented) first passage limit shapes, see \cite{durrett1981shape}, \cite{marchand2002strict}, and Section 2.5.1 of \cite{auffinger201750}.

The goal of this paper is to carry out some of the major steps of the KPZ program for the universality class observed at such transitions.
We show in Section \ref{s:SJtosBr} that along the critical line in the Sepp\"al\"ainen--Johansson model,
one observes a different universality class with exponents $0:1:2$.
The full 4-parameter scaling limit of the passage value is given by the {\bf Brownian web landscape} $\mathcal B(x,s;y,t)$, a point to half-line version of the Brownian web distance from \cite{VV26}.

\begin{figure}
\centering
\includegraphics[height=0.33\textheight]{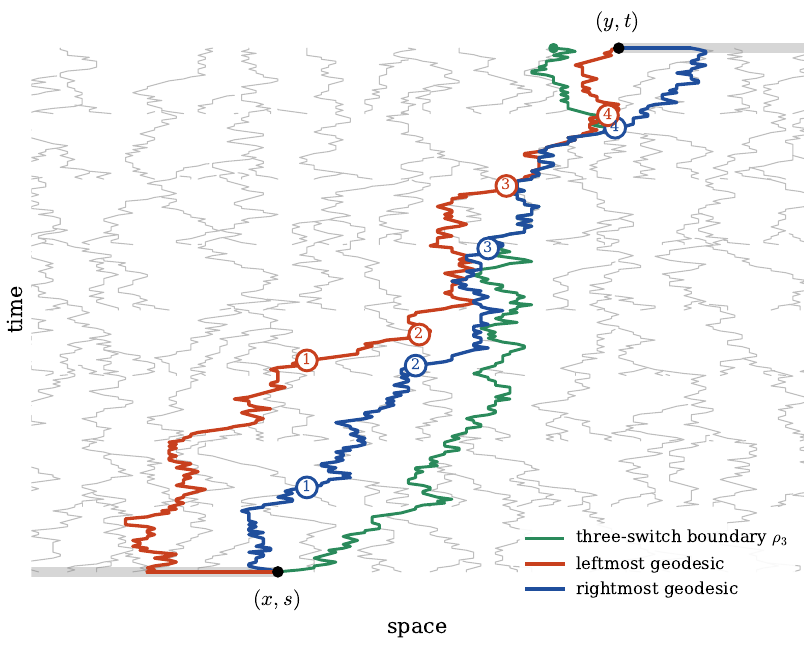}
\caption{The Brownian web landscape with $\mathcal B(x,s;y,t)=4$. The rightmost and leftmost geodesics each switch four times (numbered). The rightmost geodesic ends with a left jump to $(y,t)$, and the leftmost one starts with a left jump from $(x,s)$. The green curve $\rho_3$ bounds the points reachable from $(x,s)$ with at most three switches; it ends to the left of $y$.}
\end{figure}

We show in Section \ref{s:ips} that the evolutions of discrete-time TASEP with sequential update from left to right and from right to left as well as the geometric PushTASEP are naturally embedded into the Sepp\"al\"ainen--Johansson model.
As a consequence, our results imply limit theorems for these discrete-time versions of TASEP dynamics which can be of independent interest.

Unlike the Brownian web distance, the Brownian web landscape has a time-reversal symmetry that is related to planar duality, see Section \ref{s:duality}.
We study its large-scale behavior in Section \ref{s:large-scale}, geodesics in Section \ref{s:geodesics}, and associated height evolution in Section \ref{s:evolution}.

In Section \ref{s:large-scale}, we prove that the Brownian web landscape almost surely satisfies the large-scale limit
\begin{equation}
\frac{\mathcal B(rx,rs;ry,rt)}{r}
\to\mathcal D(x,s;y,t)=\mathbf 1_{x<y}\frac{(x-y)^2}{4(t-s)}.
\end{equation}
locally uniformly on $\{s<t\}$ as $r\to\infty$.
We also give simultaneous bounds on the deviations from $\mathcal D$
over all space-time endpoints. Up to logarithmic factors, these bounds
have size $(1+\mathcal D)^{1/3}$, with upper-tail exponent $2$
and lower-tail exponent $3/2$; see Theorem \ref{t:large-dev}.

Next, we provide the full description of geodesics of the Brownian web landscape in Section \ref{s:geodesics}. In particular, any Brownian web path, and any dual Brownian web path are geodesics. Moreover, we prove that unlike in the directed landscape, there are uncountably many geodesics between every properly ordered pair of points.

In Section \ref{s:evolution}, we define the height evolution
\begin{equation}\label{evolution}
h_t(y)=\inf_{x\in\R}\left(h_0(x)+\mathcal B(x,0;y,t)\right),
\end{equation}
the analogue of the KPZ fixed point in this class.
We establish its Markov property and
for nondecreasing lower semicontinuous integer-valued initial data
bounded from below and tending to infinity at $+\infty$,
identify its level boundaries with one-sided reflected Brownian motions.
This gives Fredholm determinant formulas for its transition probabilities through \cite{NQR20}.

We identify a family of increment-stationary laws for the height evolution \eqref{evolution}.
If $h_0$ is a lower semicontinuous two-sided Poisson counting function of intensity $\lambda>0$,
normalized by $h_0(0)=0$ and independent of the Brownian web, then the spatial increments of the height $(h_t(x)-h_t(0))_{x\in\R}$ at any time $t\ge0$ also follow the same distribution as the initial Poisson counting function.
Thus Poisson counting functions play the role of Brownian increment-stationary initial data for the KPZ fixed point.

Constructing the dynamical version of the directed landscape is still an open problem.
In the $0:1:2$ universality class, there is a natural dynamical limit object based on the dynamical Brownian web, see \cite{HW09}, and \cite{NRS10}.
In \cite{RS24} it was shown that the dynamical discrete web converges to the dynamical Brownian web.
We introduce dynamical Sepp\"al\"ainen--Johansson model in Section \ref{s:dynamics} by resampling the horizontal edge weights according to a Poisson point process in the dynamical time variable.
We use the results mentioned above on the dynamical Brownian web to show that the dynamical Sepp\"al\"ainen--Johansson model converges to the Brownian web landscape on the dynamical Brownian web in the sense of finite-dimensional distributions.

In related work, \cite{Belair25} proves convergence to the Brownian web distance for a class of discrete first passage percolation models, including models in which walks can cross before coalescing, using an argument that does not rely on planar duality.
\cite{Rahman26} studies order-one fluctuations in the transition region of the Sepp\"al\"ainen--Johansson model and obtains limiting statistics expressed by Fredholm determinants and related to the largest GUE eigenvalue.

\section{Model and main results}

\subsection{Sepp\"al\"ainen--Johansson distance}

The usual definition of the Sepp\"al\"ainen--Johansson model
is given on the two-dimensional integer lattice $\Z^2$ with edges directed up and right as follows.
Let $\xi_{ij}$ be i.i.d.\ Bernoulli random variables with parameter $1/2$ for all $i,j\in\Z$.
We assign weight $w(e)$ to each directed edge $e$ of the lattice $\Z^2$ by letting
\begin{equation}\label{edgeweights}
w(e)=\left\{\begin{array}{cl} \xi_{ij} & \mbox{if $e$ is the directed horizontal edge from $(i,j)$ to $(i+1,j)$,} \\
0 & \mbox{if $e$ is the directed vertical edge from $(i,j)$ to $(i,j+1)$} \end{array}\right.
\end{equation}
for all $i,j\in\Z$.
An admissible path $\pi:(i,j)\to(k,l)$ is a sequence of directed edges $e_1,e_2,\dots,e_{k+l-i-j}$
leading from $(i,j)$ to $(k,l)$ using only up and right steps.
The weight of this path $\pi$ is given by $w(\pi)=\sum_{m=1}^{k+l-i-j}w(e_m)$.
We define the first passage value
\begin{equation}\label{defF}
F(i,j;k,l)=\min_{\pi:(i,j)\to(k,l)} w(\pi)
\end{equation}
where the minimum is taken for all admissible paths from $(i,j)$ to $(k,l)$ and it is $\infty$ if there is no admissible path.

We define the Sepp\"al\"ainen--Johansson distance as the first passage value $F$ on the Sepp\"al\"ainen--Johansson model rotated by $45^\circ$
and scaled by $\sqrt2$ which is consistent with the directional convention followed in the literature on the directed landscape.
Let
\begin{equation}
\Z^2_\e=\{(i,j)\in\Z^2:i+j\mbox{ is even}\}
\end{equation}
be the even sublattice of $\Z^2$ and we set
\begin{equation}\label{SJF}
D^\SJ(i,j;k,l)=F\left(\frac{i+j}2,\frac{j-i}2;\frac{k+l}2,\frac{l-k}2\right)
\end{equation}
for any $(i,j),(k,l)\in\Z^2_\e$.

On the same lattice $\Z^2_\e$ we define the random walk web distance $D^\RW$ using the same set of random variables $\xi_{ij}$.
The underlying random walk web consists of coalescing random walk paths with time being the vertical direction.
From a point $(i,j)\in\Z^2_\e$, we let the random walk go to $(i+1-2\xi_{\frac{i+j}2,\frac{j-i}2},j+1)$,
that is, we assign $0$ weight to this edge and weight $1$ to the edge from $(i,j)$ to $(i-1+2\xi_{\frac{i+j}2,\frac{j-i}2},j+1)$.
The random walk web distance $D^\RW(i,j;k,l)$ is the minimal weight of a directed path from $(i,j)$ to $(k,l)$ in $\Z^2_\e$
using only edges of the form $(m,n)$ to $(m\pm1,n+1)$.
The distance is infinite if there is no such path.

Intuitively, the distance $D^\RW(i,j;k,l)$ is equal to the number of jumps needed between neighbouring random walk paths
to proceed from $(i,j)$ to $(k,l)$.
The underlying random walk started at $(i,j)$ follows the rightmost points
that can be reached from $(i,j)$ with zero Sepp\"al\"ainen--Johansson distance.
We let $D^\LF(i,j;i-1,j+1)=0$ for all $(i,j)\in\Z^2_\e$ and infinite otherwise to be the left-free directed metric on $\Z^2_\e$
which allows free jumps in the up-left direction.
The directed metric induced by another given directed metric was introduced in Definition 5.3 of \cite{DV21}
as the supremum of all directed metrics that are smaller than the given one.

\begin{proposition}\label{prop:SJinduced}
\begin{enumerate}
\item
The Sepp\"al\"ainen--Johansson distance $D^\SJ$ is the directed metric
induced by the random walk web distance $D^\RW$ and left-free directed metric $D^\LF$.
\item
For $(i,j),(k,l)\in\Z^2_\e$ with $|k-i|\le l-j$, we also have
\begin{equation}\label{eq:SJhalfline}
D^\SJ(i,j;k,l)=\min_{\substack{z\ge k\\z\equiv l\pmod2}}D^\RW(i,j;z,l).
\end{equation}
\end{enumerate}
\end{proposition}
We introduce the rescaled Sepp\"al\"ainen--Johansson distance given as
\begin{equation}\label{eq:diffusivescaling}
D^\SJ_n(x,s;y,t)=D^\SJ(n^{1/2}x,ns;n^{1/2}y,nt)
\end{equation}
for all $x,s,y,t\in\R$ where we suppressed the integer parts in the above notation on the right-hand side.
We prove the convergence of the rescaled Sepp\"al\"ainen--Johansson distance to the Brownian web landscape
where the limit is defined in terms of the Brownian web.

\subsection{Brownian web distance}

The Brownian web consists of trajectories of one-dimensional Brownian motions starting at every point of the plane.
It was first constructed in \cite{TW98}, the appropriate topology and the term Brownian web were introduced in \cite{FINR04}.

Let $\Pi_0$ denote the space of all continuous paths defined on intervals of type $[t_0,\infty)$
with the topology of uniform convergence on compact intervals.
Let $\cH$ be the space of all compact collections of paths in $\Pi_0$ with the Hausdorff metric.
The Brownian web $B$ is a random element of $\cH$ with the property that if $\gamma$ defined on $[t,\infty)$ is in $B$,
then so is $\gamma$ restricted to any interval $[s,\infty)$ with $s\ge t$.
Let $\wh\cH$ denote the space of compact collections of continuous backward paths.
The dual or backward Brownian web $\wh B$ is the almost surely unique element of $\wh\cH$ which does not cross any forward path in $B$.

A path $\gamma$ in $B$ with $\gamma(t)=x$ is called an outgoing path from $(x,t)$ if it is defined on $[t,\infty)$
and the part $\gamma_{[s,t]}$ is an incoming path at $(x,t)$ if $s<t$.
Almost surely, almost all $(x,t)\in\R^2$ have a unique outgoing path which we denote by $B_{(x,t)}$.
Exceptional points of type $(1,2)$ have one incoming and two outgoing paths where
exactly one of the outgoing paths is the continuation of the incoming path.
This outgoing path is the Brownian web path $B_{(x,t)}$ starting from $(x,t)$ and the additional path is denoted by $B'_{(x,t)}$.
For other types of exceptional points with multiple outgoing paths, a unique $B_{(x,t)}$ is not defined.
We denote by $B_{(x,t)+}$ and $B_{(x,t)-}$ the rightmost and leftmost outgoing paths from $(x,t)$, that is, the supremum and infimum of all outgoing paths from $(x,t)$ in $B$, respectively.
For the dual web, $\wh B_{(x,t)+}$ and $\wh B_{(x,t)-}$ denote the supremum and infimum of the backward outgoing paths starting at $(x,t)$ in $\wh B$, respectively.

The Brownian web distance is an integer-valued scale-invariant random directed metric on the plane introduced in \cite{VV26}
defined in terms of the Brownian web.
In this paper, we only use the lower semicontinuous version of the Brownian web distance denoted by $D^{\Br,\LSC}$ in \cite{VV26}.
For simplicity, we use the notation $D^\Br$ for the lower semicontinuous Brownian web distance in this paper which we define as follows.

The lower semicontinuous Brownian web distance $D^\Br(x,s;y,t)$ is the smallest integer $k\ge 0$
for which there exist points $(x,s)=(x_0,t_0),(x_1,t_1),\dots,(x_{k+1},t_{k+1})=(y,t)\in\R^2$ with $t_0<t_1<t_2<\dots<t_{k+1}$
and a continuous path $\pi:[s,t]\to\R$ with $\pi(t_i)=x_i$ for $i=0,\dots,k+1$ so that for each $i=0,\dots,k$ there is a path $\gamma_i$ in the Brownian web
with $\pi(r)=\gamma_i(r)$ for all $r\in [t_i,t_{i+1}]$.
We set $D^\Br(x,s;y,t)=+\infty$ if there is no such $k\ge 0$.

We denote by $\mathcal B$ the \emph{Brownian web landscape} defined as
\begin{equation}\label{defDsBr}
\mathcal B(x,s;y,t)=D^\Br(x,s;[y,\infty),t)=\inf_{z\ge y}D^\Br(x,s;z,t),\qquad s<t.
\end{equation}
We set $\mathcal B(x,s;y,s)=0$ if $y\le x$, and $\mathcal B(x,s;y,t)=\infty$ if $s>t$ or if $s=t$ and $y>x$.
The version used for metric composition is introduced in Section~\ref{s:evolution}.

We mention that one can introduce a left (or a right) version of the Brownian web distance
which only allows jumps to the left (or right) at those type $(1,2)$ point
where the additional outgoing path is to the left (or right) of the continuation of the incoming path.

The Brownian web landscape $\mathcal B$ preserves the symmetry of the Sepp\"al\"ainen--Johansson model
under the rotation by $180^\circ$, that is,
\begin{equation}
F(i,j;k,l)\stackrel\d=F(-k,-l;-i,-j)
\end{equation}
holds for all $i,j,k,l\in\Z$ simultaneously.
The dual symmetry of $\mathcal B$ below is proved in Section \ref{s:duality}.
It uses the Brownian web distance $\wh D^\Br$ defined on the dual Brownian web $\wh B$.
The dual Brownian web after time reversal has the same law as the forward web
and the dual Brownian web distance $\wh D^\Br(y,t;x,s)$ can similarly be defined as the minimal number of jumps at $(1,2)$ points
to get from $(y,t)$ to $(x,s)$ using paths of the dual Brownian web $\wh B$.

\begin{proposition}\label{prop:duality}
Jointly for all $(x,s,y,t)\in\R^4$ with $s<t$ we have that
\begin{equation}\label{duality}
\mathcal B(x,s;y,t)=\wh D^\Br(y,t;(-\infty,x],s)\stackrel{\d}=\mathcal B(-y,-t;-x,-s)
\end{equation}
where the first equality holds almost surely and the second one in distribution.
The distributional symmetry extends to $s\ge t$ by the definition of $\mathcal B$.
\end{proposition}

\subsection{Discrete-time TASEP}

We establish a connection between the Sepp\"al\"ainen--Johansson model and two versions of the discrete-time TASEP with sequential update as well as the geometric PushTASEP.
The discrete-time totally asymmetric simple exclusion process (TASEP) with sequential update was first introduced in \cite{RS05},
and its asymptotic fluctuations were identified in \cite{BFP07}.
Various update rules for discrete-time TASEP were discussed in \cite{MS11}.
In versions of discrete-time TASEP, the order of updates is crucial.
The order from right to left corresponds to rule 1 in \cite{MS11}, the update from left to right corresponds to rule 2.
These two rules have a natural symmetry obtained by exchanging left and right and exchanging the roles of particle and hole.
Moreover, the only translation-invariant stationary ergodic measures with constant marginals for these dynamics
are the Bernoulli product measures by Proposition 3.1 in \cite{MS11}.
The geometric PushTASEP was introduced in \cite{MR23} where explicit formulas are given for the particle positions.

The discrete-time TASEP with sequential update from right to left on $\Z$ is described by the particle locations $(x_j(i))_{i,j}$
where $x_j(i)$ denotes the position of the $j$th particle at time $i$ for all $j=0,1,2,\dots$ and $i=0,1,2,\dots$.
Particle configurations allow for at most one particle at each integer site
and we use the convention that particles are numbered from right to left, that is, $x_j(i)<x_{j-1}(i)$.
Given the particle locations at time $i$, the new locations at time $i+1$
are updated sequentially with $j=0,1,2,\dots$ increasingly as follows.
Each particle attempts to jump to the right by one with probability $1/2$ with all attempts independent for all particles and all times.
The jump attempt of the $j$th particle is successful, that is, $x_j(i+1)=x_j(i)+1$ if the target location is empty, that is,
if $x_j(i)+1<x_{j-1}(i+1)$ and it stays put otherwise with the convention $x_{-1}(i)=\infty$ for all $i$.
The ordering of the particles is preserved by the dynamics.
Step initial condition means $x_j(0)=-j$ for $j=0,1,2,\dots$.

The discrete-time TASEP with sequential update from left to right is defined similarly.
Let $y_j(i)$ denote the position of the $j$th particle at time $i$ for all $j=0,1,2,\dots$ and $i=0,1,2,\dots$.
There is at most one particle at each location and the particles are numbered from right to left, $y_j(i)<y_{j-1}(i)$.
Given the particle locations at time $i$, the new locations at time $i+1$
are updated sequentially with $j=0,1,2,\dots$ decreasingly as follows.
The $j$th particle with a gap of size $m$ in front of it jumps to the right by the minimum of $m$ and a geometric random variable of parameter $1/2$,
that is, it jumps by $l$ with probability $2^{-l-1}$ for $l\in\{0,\dots,m-1\}$ and by $m$ with probability $2^{-m}$.
We assume that all jumps are independent for all particles and all times.
The update dynamics is well-defined for finite configurations and for step initial condition $y_j(0)=-j$ with $j=0,1,2,\dots$,
moreover, it preserves the ordering of the particles.

In geometric PushTASEP on $\Z$ we denote by $w_i(j)$ the position of the $i$th particle at time $j$ for all $i=0,1,2,\dots$ and $j=0,1,2,\dots$.
There is always at most one particle at each position and particles are numbered from left to right meaning that $w_{i-1}(j)<w_i(j)$.
Particle positions are updated from left to right as follows.
Given the locations at time $j-1$ for all particles and at time $j$ for particles $0,1,\dots,i-1$,
the $i$th particle is pushed to the right of the $(i-1)$st particle if needed and then it jumps to the right by a geometric amount.
Formally, positions are updated using i.i.d.\ geometric variables $G_{ij}$ with distribution $P(G_{ij}=k)=2^{-k-1}$ for $k=0,1,2,\dots$
and we let $w_i(j)=\max(w_i(j-1),w_{i-1}(j)+1)+G_{ij}$.

In the particle system studied in \cite{S98} the location of particle $i$ at time $j$ is given by $z_i(j)$ for $i\in\Z$ and $j=0,1,\dots$.
Particle updates are defined in (3.4) of \cite{S98} in terms of a last passage value given by the background variables of the Sepp\"al\"ainen--Johansson model that we denote by $\xi_{ij}$ in this paper.
The updates are proved to satisfy the semigroup property, hence configurations evolve as a Markov chain.
It follows from the definition of updates that if all sites of $\Z_-$ are initially occupied, that is, $z_i(0)=i$ for all $i\le0$ then these particles stay put forever.

The evolution of the system $(z_i(j))_{i,j}$ is described in (3.7) of \cite{S98} explicitly as follows.
Given the positions $z_i(j-1)$ and $z_{i-1}(j)$, we have that
\begin{equation}\label{zijevo}
z_i(j)=\min(z_i(j-1),b_{i-1}(j))
\end{equation}
where
\begin{equation}
b_{i-1}(j)=\min\{m>z_{i-1}(j):\xi_{m-1,j}=0\}
\end{equation}
is the first position after the particle $z_{i-1}(j)$ which is the right endpoint of a free horizontal edge.
By the property that the infinite block of particles on $\Z_-$ does not move, it is enough to keep track of particles of nonnegative index, that is $z_i(j)$ for $i=0,1,\dots$ and $j=0,1,\dots$, and the evolution rule \eqref{zijevo} can be applied in $i$ increasingly.

\subsection{Dynamical distances}

Let $Y$ be an instance of the discrete web, that is, we have a collection of i.i.d.\ Bernoulli random variables $(\xi_{i,j})_{i,j \in \Z}$
which define the edge weights in the Sepp\"al\"ainen--Johansson model in \eqref{edgeweights},
and for each lattice point $(x,t)\in\Z^2_\e$, we draw an edge going from $(x,t)$ to $(x+1-2\xi_{\frac{t+x}2,\frac{t-x}2},t+1)$.
Then $Y$, as the collection of all random walk trajectories starting from every point of $\Z^2_\e$, is considered as an $\Hc$-valued random variable.
Let $(N_{x,t}(\tau))_{(x,t) \in \Z^2_\e}$ be a collection of i.i.d.\ Poisson processes of intensity $\lambda>0$ on the positive real line.
At each point in the process $N_{x,t}(\tau)$, we switch the direction of the arrow starting at $(x,t)$.
Setting $Y_0=Y$, we obtain a process $(Y(\tau))_{\tau \ge 0}$ taking its values in $\Hc$.
This object is called the \emph{dynamical discrete web} (DyDW).

The DyDW was first introduced by Howitt and Warren in \cite{HW09}.
They conjectured that under diffusive scaling, by taking the intensity of the switch rates to depend on $n$ via $\lambda=1/\sqrt{n}$,
then there is a non-trivial scaling limit.
Assuming the existence of this limit, they also gave an explicit description of its finite-dimensional distributions in terms of sticky Brownian motions.
A process with these finite-dimensional distributions was constructed in \cite{NRS10} using a special Poissonian marking of the $(1,2)$ points.
This process is now known as the \emph{dynamical Brownian web} (DyBW).
Convergence of the DyDW to the DyBW as a process was shown more recently in \cite{RS24}.
More precisely, they showed that there is a version of the DyBW which has c\`adl\`ag trajectories when viewed as a random function from $[0,\infty)$ to $\Hc$,
and that the appropriately rescaled DyDW converges in distribution to the DyBW in the Skorokhod topology on c\`adl\`ag functions.

We define dynamical versions of the Sepp\"al\"ainen--Johansson, left-free and Brownian web distances,
$D^{\bullet}(x,s;y,t;\tau)$, and of the Brownian web landscape, $\mathcal B(x,s;y,t;\tau)$, by evaluating the corresponding function at dynamical time $\tau$ in the DyDW or DyBW.
These are all stationary in $\tau$, since the DyDW and DyBW are trivially so.
The interest in studying these dynamical distances comes from the fact that there are random exceptional times
when the DyDW and DyBW have atypical behaviour from the usual discrete web and Brownian web.
We expect this unusual behaviour at random times to manifest itself in the dynamical distances.

We use base switching rate $\lambda=1$ in the unscaled dynamical model.
Similarly to \eqref{eq:diffusivescaling}, we introduce the diffusively rescaled  dynamical random walk web distance and Sepp\"al\"ainen--Johansson distance as
\begin{align}
D_n^\RW(x,s;y,t;\tau)&=D^\RW(n^{1/2}x,ns;n^{1/2}y,nt;n^{-1/2}\tau),\label{defDRWntau}\\
D_n^\SJ(x,s;y,t;\tau)&=D^\SJ(n^{1/2}x,ns;n^{1/2}y,nt;n^{-1/2}\tau)\label{defDSJntau}
\end{align}
and for any $\tau\ge0$ we denote by $D_n^{\RW}(\tau)$ and $D_n^{\SJ}(\tau)$ the distance functions above
for all pairs of points $(x,s),(y,t)\in\R^2$ at dynamical time $\tau$.

\subsection{Results}

First we show that in the appropriate coupling, the rescaled static Sepp\"al\"ainen--Johansson model converges
to the Brownian web landscape almost surely in the epigraph sense.
We define the epigraph of the lower semicontinuous function $f:\R^4\to\ol\R=\R\cup\{\infty\}$ as the closed set
\begin{equation}
\epi f=\{(x,y)\in\R^4\times\ol\R:y\ge f(x)\}.
\end{equation}
The set $\mathcal E_*$ of all closed subsets $\Gamma\subset\R^4\times\ol\R$ such that $\Gamma\cap(\{x\}\times\ol\R)\neq\emptyset$ for all $x\in\R^4$
contains the epigraphs of all lower semicontinuous functions.
The space $\mathcal E_*$ is equipped with a version of the Hausdorff topology.
Let $E:\R^4\times\ol\R\to\R^4\times[-1,1]$ be the map given by
\begin{equation}
E(r,s)=\left(r,\frac{|u-v|e^{-|r|}s}{1+|s|}\right)
\end{equation}
for any $(r,s)\in\R^4\times\ol\R$ with $r=(u,v)\in\R^4$ with the convention $E(r,\pm\infty)=(r,\pm|u-v|e^{-|r|})$.
For any two elements $e_1,e_2\in\mathcal E_*$,
we define their distance $d_*(e_1,e_2)$ as the Hausdorff distance of the images $E(e_1)$ and $E(e_2)$.
The space $(\mathcal E_*,d_*)$ is compact by Lemma 7.1 of~\cite{DV21}.

\begin{theorem}\label{thm:SJtosBr}
There is a coupling of the Brownian web and the sequence of Sepp\"al\"ainen--Johansson models
so that the epigraph of the lower semicontinuous function $D^\SJ_n$ on $\R^4$ converges to the epigraph of $\mathcal B$,
that is, $\epi D^\SJ_n\to\epi \mathcal B$ in $\mathcal E_*$ almost surely as $n\to\infty$.
\end{theorem}

The first passage values $F(0,0;i,j)$ from the origin to targets in $\Z_+^2$ in the Sepp\"al\"ainen--Johansson model
can naturally be described in terms of interacting particle systems.
We provide four different particle systems embedded in the model.
The first two systems described in Proposition \ref{prop:dTASEP} are the discrete-time TASEP with sequential update from right to left and from left to right.
They arise when time is described by index $i$, that is, a configuration of the systems corresponds to the passage times in the Sepp\"al\"ainen--Johansson model along a vertical line.
The other two systems given in Proposition \ref{prop:SeppIPS} are the geometric PushTASEP and the particle system considered in \cite{S98}.
In this case, time is index $j$, that is, time slices are horizontal.

\begin{proposition}\label{prop:dTASEP}
\begin{enumerate}
\item
With the location of the $j$th particle at time $i$ given by
\begin{equation}\label{defxij}
x_j(i)=F(0,0;i,j)-j
\end{equation}
for all $i=0,1,2,\dots$ and $j=0,1,2,\dots$,
the particle system evolves as discrete-time TASEP with sequential update from right to left
started from step initial condition $x_j(0)=-j$.
The indicator of the jump attempt of the $j$th particle at time $i$ is $\xi_{ij}$.

\item
For all $i=0,1,2,\dots$ and $j=0,1,2,\dots$, we define
\begin{equation}\label{defyij}
y_j(i)=F^{-1}(0,0;i,j)-j=\min\{k\in\Z:F(0,0;i,k)\le j\}-j
\end{equation}
as the position of the $j$th particle at time $i$.
Then the particle system $y_j(i)$ evolves as discrete-time TASEP with sequential update from left to right
started from step initial condition $y_j(0)=-j$.
The $j$th particle jumps to the right by $l$ at time $i$ if and only if $l$ is the largest integer in
$\{0,1,\dots,y_{j-1}(i)-y_j(i)-1\}$ with $\xi_{i,y_j(i)+j}=\dots=\xi_{i,y_j(i)+j+l-1}=1$.

\item
For any time $i=0,1,\dots$, the locations of the holes between particles at $(x_j(i))_{j=0}^\infty$ are given by the positions $(1-y_j(i))_{j=0}^\infty$.
\end{enumerate}
\end{proposition}

By Proposition \ref{prop:dTASEP}, the convergence result for the Sepp\"al\"ainen--Johansson model can directly be translated to a limit theorem for
discrete-time TASEP with sequential update from right to left on the left edge of its rarefaction fan
and for particle flux for left to right update.
We mention that a limit theorem for geometric PushTASEP and for the system $(z_i(j))_{i,j}$ can be deduced similarly.

\begin{corollary}\label{cor:dTASEPconv}
\begin{enumerate}
\item 
For the distance of particles from their initial position in discrete-time TASEP with right to left sequential update, we have
\begin{equation}\label{dTASEPconv}
x_{\lfloor n-\sqrt2y\sqrt n\rfloor}(n)-x_{\lfloor n-\sqrt2y\sqrt n\rfloor}(0)\stackrel\d\Longrightarrow \mathcal B(0,0;y,1)
\end{equation}
as $n\to\infty$ for each fixed $y\in\R$.
In particular,
\begin{equation}\label{dTASEP0}
\lim_{n\to\infty}P\left(x_{\lfloor n-\sqrt2y\sqrt n\rfloor}(n)=x_{\lfloor n-\sqrt2y\sqrt n\rfloor}(0)\right)=\int_y^\infty\frac1{\sqrt{2\pi}}e^{-z^2/2}\,\d z.
\end{equation}

\item
Let
\begin{equation}
N(i,k)=|\{j=0,1,2,\dots:y_j(i)-y_j(0)\ge k\}|=|\{j=0,1,2,\dots:y_j(i)\ge k-j\}|
\end{equation}
denote the number of particles in the discrete-time TASEP with left to right sequential update that have made at least $k$ steps by time $i$.
Then
\begin{equation}\label{dTASEPflux}
N\left(n,\lfloor n-\sqrt2y\sqrt n\rfloor\right)\stackrel\d\Longrightarrow \mathcal B(0,0;y,1)
\end{equation}
as $n\to\infty$ for each fixed $y\in\R$.
\end{enumerate}
\end{corollary}

Since \eqref{dTASEP0} changes continuously from $1$ to $0$ as $y$ ranges over $\R$,
the left edge of the rarefaction fan after $n$ steps can be found in the window given by particles of label $n+\O(\sqrt n)$.
Corollary \ref{cor:dTASEPconv} above is therefore a limit theorem for the number of jumps on the left edge of the rarefaction fan,
which could be of independent interest.

\begin{proposition}\label{prop:SeppIPS}
\begin{enumerate}
\item
Let
\begin{equation}
w_i(j)=F^{-1}_3(0,0;i,j)=\max\{k\in\Z:F(0,0;k,j)\le i\}\label{defwij}
\end{equation}
denote the position of the $i$th particle at time $j$ for $i=0,1,2,\dots$ and $j=0,1,2,\dots$.
Then the system follows the evolution of geometric PushTASEP and we have the equality of events
\begin{equation}\label{wevolution}
\{w_i(j)\ge M+m\}=\{\xi_{M,j}=\dots=\xi_{M+m-1,j}=0\}
\end{equation}
for all $m=1,2,\dots$ where
\begin{equation}\label{defM}
M=\max(w_i(j-1),w_{i-1}(j)+1).
\end{equation}
We use the conventions $w_i(-1)=0$ for $i\ge0$ and $w_{-1}(j)=-1$ for $j\ge0$.

\item
We let
\begin{equation}\label{defzij}
z_i(j)=\min\{z\in\Z:z-F(0,0;z,j)\ge i\}
\end{equation}
the position of the $i$th particle at time $j$ for $i=0,1,2,\dots$ and $j=0,1,2,\dots$ and we let $z_i(j)=i$ for all $i\le0$ and all $j\ge0$.
Then the system evolves as the particle system in \cite{S98}, that is, given by \eqref{zijevo}.
The initial condition is given by $z_i(0)=i$ for $i\le0$ and $z_i(0)=p_i$ for $i\ge1$ where $p_i$ is the $i$th position $m\ge1$ with $\xi_{m-1,0}=0$.

\item
The holes between particles at $(z_i(j))_{i,j}$ are located at $(1+w_i(j))_{i=0}^\infty$ for all time $j$.
In particular, the evolution of the holes in the particle system in \cite{S98} follows the evolution of PushTASEP.
\end{enumerate}
\end{proposition}

Our results also describe the large-scale behavior, height evolution,
and geodesics of the Brownian web landscape.
For $s<t$, let
\begin{equation}
\mathcal D(x,s;y,t)=\mathbf 1_{x<y}\frac{(x-y)^2}{4(t-s)}
\end{equation}
be the one-sided Dirichlet distance.
The following theorem gives simultaneous bounds on the deviations
of $\mathcal B$ from $\mathcal D$.

\begin{theorem}\label{t:large-dev}
There exist random constants $C_1,C_2\ge 0$
and a deterministic constant $c>0$ such that,
almost surely, for every $u=(x,s;y,t)\in\R^4$ with
$s<t$,
\begin{equation}
-C_1(1+\D(u))^{1/3}\ell(u)^{2/3}
\le (\mathcal B-\D)(u)
\le C_2(1+\D(u))^{1/3}\ell(u)^{1/2},
\end{equation}
where
\begin{equation}
\ell(u)=\log\left(2+\|u\|_\infty+\frac1{t-s}\right).
\end{equation}
Moreover, for every $a\ge0$,
\begin{equation}
P(C_1\ge a)\le e^{1-ca^{3/2}},
\qquad
P(C_2\ge a)\le e^{1-ca^2}.
\end{equation}
\end{theorem}

As a consequence, we obtain the following locally uniform law of
large numbers.

\begin{corollary}\label{cor:LLN}
Almost surely, for every compact set
$K\subset\{(x,s;y,t)\in\R^4:s<t\}$,
\begin{equation}
\lim_{r\to\infty}
\sup_{(x,s;y,t)\in K}
\left|
\frac{\mathcal B(rx,rs;ry,rt)}{r}
-\D(x,s;y,t)
\right|=0.
\end{equation}
\end{corollary}

Theorem \ref{t:large-dev} and Corollary \ref{cor:LLN} are proved
in Section \ref{s:large-scale}.

The Brownian web landscape defines a height evolution by \eqref{evolution}.
Poisson counting functions are increment-stationary for this evolution,
just as Brownian motion is for the KPZ fixed point.

\begin{theorem}\label{prop:poissonstationarity}
Let $\lambda>0$, and let $h_0$ be a lower semicontinuous
two-sided Poisson counting function of intensity $\lambda$,
normalized by $h_0(0)=0$ and independent of the Brownian web.
Then, for every $t\ge0$,
\begin{equation}
\bigl(h_t(x)-h_t(0)\bigr)_{x\in\R}
\stackrel{\d}=
\bigl(h_0(x)\bigr)_{x\in\R}.
\end{equation}
\end{theorem}

We prove Theorem \ref{prop:poissonstationarity} in
Section \ref{s:evolution}, where we also establish the Markov property
of the height evolution and give Fredholm determinant formulas
for its transition probabilities.

Next we describe the geodesics of the Brownian web landscape.
For a more precise statement, see Theorem \ref{thm:geodesicchar}.

\begin{theorem}\label{thm:geod_intro}
Let $(x,s,y,t)\in\R^4$ with $s<t$ and $\mathcal B(x,s;y,t)=k$.
The rightmost geodesic $\pi$ and the leftmost geodesic $\wh\pi$ between $(x,s)$ and $(y,t)$
can be described explicitly in terms of the forward and backward Brownian webs.

Let $\gamma$ be a c\`adl\`ag path between $\wh\pi$ and $\pi$ which only has jumps down.
It is a geodesic between $(x,s)$ and $(y,t)$ if and only if
the modulus of continuity of $\gamma$ can locally be upper bounded by that of the Brownian web except for $k$ points
and by that of the rightmost outgoing Brownian web path at the exceptional points.
The $k$ exceptional points are exactly the points of increase of the geodesic length.
\end{theorem}

The rightmost and leftmost geodesics are two different paths almost surely and there are uncountably many geodesics between them.

\begin{corollary}\label{cor:infgeod}
Almost surely, for all $(x,s,y,t)\in\R^4$ with $s<t$, there are uncountably many geodesics between $(x,s)$ and $(y,t)$.
\end{corollary}

The next result provides a generalization of Theorem \ref{thm:SJtosBr} for the dynamical setting.

\begin{theorem}\label{thm:fdconvdyn}
Fix the deterministic dynamical times $0\le\tau_1<\dots<\tau_k$.
There exists a coupling of the dynamical Sepp\"al\"ainen--Johansson models and the Brownian web at these dynamical times such that
\begin{equation}\label{fdconvRWtoBr}
(D_n^{\RW}(\tau_1),\dots,D_n^{\RW}(\tau_k))\to(D^\Br(\tau_1),\dots,D^\Br(\tau_k))
\end{equation}
and
\begin{equation}\label{fdconvSJtosBr}
(D_n^{\SJ}(\tau_1),\dots,D_n^{\SJ}(\tau_k))\to(\mathcal B(\tau_1),\dots,\mathcal B(\tau_k))
\end{equation}
almost surely in the product topology on $(\mathcal{E}_*)^k$
where $(D_n^{\RW}(\tau))_{\tau \ge 0}$ and $(D_n^{\SJ}(\tau))_{\tau \ge 0}$ are the rescaled dynamical random walk web distances
and Sepp\"al\"ainen--Johansson distances given in \eqref{defDRWntau}--\eqref{defDSJntau}.
\end{theorem}

This theorem gives convergence of the finite-dimensional distributions of the dynamical random walk web distance and Sepp\"al\"ainen--Johansson distance
to the dynamical Brownian web distance and Brownian web landscape, respectively.
We conjecture that these convergences also hold as processes in the same sense as in \cite{RS24},
that is, in the Skorokhod topology on c\`adl\`ag functions from $[0,\infty)$ to $\mathcal{E}_*$.

\section{Convergence to the Brownian web landscape}\label{s:SJtosBr}

We prove Proposition \ref{prop:SJinduced} and Theorem~\ref{thm:SJtosBr} in this section.

\begin{proof}[Proof of Proposition \ref{prop:SJinduced}]
\begin{enumerate}
\item
The minimum of the random walk web edge weights and the left-free edge weights is exactly the SJ edge weight, proving the induced-metric assertion.
\item
We first prove that the left-hand side of \eqref{eq:SJhalfline} is an upper bound for the right-hand side.
Given an SJ path $\pi$ from $(i,j)$ to $(k,l)$, we construct a walk-web path weakly to its right as follows.
When strictly to its right, we follow the zero-cost arrow.
At a meeting, we follow that arrow unless $\pi$ steps right and the arrow points left; in this case we take the right edge at cost one.
Each charge is also present in $\pi$, so the constructed path ends at some $z\ge k$ with no greater cost.

Conversely, given a walk-web path $\pi'$ from $(i,j)$ to $(z,l)$ with $z\ge k$, we show that there is a path with SJ cost at most the walk-web cost of $\pi'$.
Since $i\le k+l-j$ by assumption, $\pi'$ intersects the lattice line $u\mapsto k+l-u$.
For the new path, we follow $\pi'$ to the first intersection with the line $u\mapsto k+l-u$, and then we take only up-left edges to $(k,l)$.
The final edges have zero SJ cost, and the initial segment has SJ cost at most the walk-web cost of $\pi'$.
This proves \eqref{eq:SJhalfline}.
\end{enumerate}
\end{proof}

Theorem \ref{thm:SJtosBr} is a deterministic consequence of the epigraph convergence of the random walk web distance
to the Brownian web distance proved in~\cite{VV26}.
It follows from Theorem~\ref{thm:convtransform} which is a more general result
showing that the epigraph convergence of distances implies the epigraph convergence of the half-line versions.

Let $\cH$ be the space of compact collections of paths, and $\dH$ be the Hausdorff distance on $\cH$.
Thus the Brownian web and the discrete random walk web are $\cH$-valued random variables.
Given any $C \in \cH$, we can define a directed metric on $\R^4$ in the same way as the Brownian web distance is defined.
For two points $(x,s)$ and $(y,t)$ with $s<t$, their distance is the smallest integer $k$
such that there are paths $f_1,\dots,f_{k+1}$ in $C$ and times $s<r_1<\dots<r_k<t$
such that $f_1(s)=x$, $f_{k+1}(t)=y$ and $f_i(r_i) = f_{i+1}(r_i)$ for all $i$.
If no such $k$ exists, their distance is defined to be infinite.

The distance above has \emph{misalignment level $\varepsilon$} if the condition $f_i(r_i)=f_{i+1}(r_i)$ is relaxed to
\begin{equation}
|f_{i+1}(r_i)-f_i(r_i)|\le\varepsilon
\end{equation}
for all $i=1,\dots,k$.
The Brownian web distance defined with misalignment level $0$, the rescaled random walk web distance and the rescaled SJ distance given in \eqref{eq:diffusivescaling} have misalignment level $2/\sqrt n$.
With this definition, paths do not use edges with random walk web cost one but they continue on the free edge and then jump to a neigbouring site of the lattice causing a misalignment.

Given the distance $D$ defined as above on $C$, we can also define the corresponding point to half-line distance as
\begin{equation}\label{defDtilde}
\wt D(x,s;y,t)=\min_{z\geq y}D(x,s;z,t).
\end{equation}
The minimum is always achieved at some point $z=y^*$ because $D$ is nonnegative integer (or infinite) valued.
The next theorem on the open set $\{s<t\}$ provides the essential deterministic input for the convergence in Theorem \ref{thm:SJtosBr}.
The equal-time boundary for the SJ model can be treated trivially based on the following conditions for the epigraph convergence.
The sufficiency of the conditions given below are proved in Lemma 4.1 of \cite{VV26}, the necessity is a direct consequence of the definition.

\begin{lemma}\label{lemma:epigraphconv}
Let $f,f_n:\R^4\to\ol\R$ be lower semicontinuous functions.
The convergence $f_n\to f$ holds in the epigraph sense, that is, $\epi f_n\to\epi f$ in $\mathcal E_*$ as $n\to\infty$ if and only if the following two conditions are satisfied.

First, for any convergent sequence $(x_n,s_n,y_n,t_n)\to(x,s,y,t)$ it holds that
\begin{equation}\label{epiconv1}
\liminf_{n\to\infty}f_n(x_n,s_n,y_n,t_n)\ge f(x,s,y,t).
\end{equation}
Second, for any $(x,s,y,t)\in\R^4$ there exists a convergent sequence $(x_n,s_n,y_n,t_n)\to(x,s,y,t)$ such that
\begin{equation}\label{epiconv2}
\limsup_{n\to\infty}f_n(x_n,s_n,y_n,t_n)\le f(x,s,y,t).
\end{equation}
\end{lemma}

\begin{theorem}\label{thm:convtransform}
Let $C_n$ be a sequence in $\cH$ which converges to $C$ in the Hausdorff metric.
Let $D_n$ and $D$ be a nonnegative integer-valued (or infinite) distance function defined on $C_n$ and $C$ with misalignment level $\varepsilon_n$ and $0$ respectively.
Suppose that $D_n$ converges to $D$ in the epigraph topology on $\{s<t\}$ and that $\varepsilon_n\to0$.
Then the corresponding point to half-line distances $\wt D_n$ converge in epigraph to $\wt D$ on $\{s<t\}$.
\end{theorem}

For the proof of Theorem~\ref{thm:convtransform}, we state two lemmas about the Hausdorff metric and uniform convergence.

\begin{lemma}\label{lemma:sequenceHdorff}
Let $(S,d)$ be a metric space, and let $C_n$ be a sequence of compact sets in $S$ which converges in Hausdorff distance to a compact set $C$.
Suppose $x_n$ is a sequence such that $x_n\in C_n$ for each $n$.
Then there is an $x\in C$ and a subsequence $x_{n_j}$ such that $x_{n_j}\to x$ in $d$.
\end{lemma}

\begin{proof}[Proof of Lemma~\ref{lemma:sequenceHdorff}]
Let $\delta_n = \dH(C_n, C)$.
Then by definition of the Hausdorff distance, we have $C_n \subset C^{\delta_n + 1/n}$ where for a set $A \subset S$ we write
\begin{equation}
A^\delta = \{ x \in S : \text{there is } y \in A \text{ with } d(x, y) < \delta \}.
\end{equation}
So there exists a sequence $y_n$ in $C$ such that for all $n$, $d(x_n, y_n) < \delta_n + 1/n$. Since $C$ is compact, $y_n$ has a convergent subsequence in $C$, say $y_{n_j} \to x$. Then
\begin{equation}
d(x_{n_j}, x) \leq d(x_{n_j}, y_{n_j}) + d(y_{n_j}, x) \to 0.
\end{equation}
\end{proof}

\begin{proof}[Proof of Theorem~\ref{thm:convtransform}]
For a fixed $(x,s,y,t)\in\R^4$ with $s<t$, we verify the two conditions in Lemma \ref{lemma:epigraphconv}.
Property \eqref{epiconv2} is easier to check in this case, so we start with this.
By the fact that $D$ is integer-valued, the minimum is achieved in the definition \eqref{defDtilde} of $\wt D$, i.e., there is a $y^*\geq y$ such that $\wt D(x,s;y,t)=D(x,s;y^*,t)$.
Since we assume that $D_n\to D$ in epigraph, then property \eqref{epiconv2} holds for $D_n$: there exists a sequence $(x_n,s_n;y_n^*,t_n)\to(x,s;y^*,t)$ such that
\begin{equation}
\limsup_{n\to\infty}D_n(x_n,s_n;y_n^*,t_n)\leq D(x,s;y^*,t).
\end{equation}
Let $y_n = y_n^* - y^* + y$.
Then $(x_n,s_n;y_n,t_n)\to(x,s;y,t)$, and $y_n^*\geq y_n$, so
\begin{equation}
\limsup_{n\to\infty}\wt D_n(x_n,s_n;y_n,t_n)\leq\limsup_{n\to\infty}D_n(x_n,s_n;y_n^*,t_n)\leq D(x,s;y^*,t) =\wt D(x,s;y,t).
\end{equation}

Next we show property \eqref{epiconv1}.
Let $(x_n,s_n;y_n,t_n)\to(x,s;y,t)$, and set
\begin{equation}\label{LliminfD}
L = \liminf_{n\to\infty}\wt D_n(x_n,s_n;y_n,t_n).
\end{equation}
We can assume that $L$ is finite, otherwise there is nothing to show.
Since $\wt D$ and $\wt D_n$ are non-negative integer (or infinite) valued, this means that we can find a subsequence $n_j$
such that $\wt D_n(x_n,s_n;y_n,t_n)=L$ for every $n$ along the subsequence.
We redefine the sequence by taking only its elements along this subsequence and we omit subscripts.

By the path representation of \eqref{LliminfD}, there are at most $L+1$ web-path pieces joining $(x_n,s_n)$ to $(y_n^*,t_n)$, where $y_n^*\ge y_n$, with misalignment at most $\varepsilon_n$ between them.
After repeating pieces if necessary, denote their web paths by $f_1^{(n)},\dots,f_{L+1}^{(n)}$, and the times of the $i$th connection by $r_i^{(n)}$.
Then $s_n\le r_1^{(n)}\le\dots\le r_L^{(n)}\le t_n$ and
\begin{equation}
\bigl|f_i^{(n)}(r_i^{(n)})-f_{i+1}^{(n)}(r_i^{(n)})\bigr|\le\eps_n
\end{equation}
for all $i=1,\dots,L$.
Here $f_1^{(n)}(s_n)=x_n$ and $f_{L+1}^{(n)}(t_n)=y_n^*$.
By Lemma~\ref{lemma:sequenceHdorff} and compactness, we can pass to a subsequence on which the finitely many web paths and all connection times converge: $f_i^{(n)}\to f_i$ and $r_i^{(n)}\to r_i$.
Paths are extended constantly before their starting times when comparing them on a common time interval.
Successively from the first piece, whose initial value tends to $x$, convergence gives a finite limiting endpoint; $\varepsilon_n\to0$ gives the same limiting initial value for the next piece.
Thus all these limits are finite on the relevant compact time interval, and the convergence is uniform.
We can delete pieces whose time intervals collapse.
Uniform convergence and continuity give
\begin{equation}
f_i(r_i)-f_{i+1}(r_i)=\lim_{n\to\infty}\bigl(f_i^{(n)}(r_i^{(n)})-f_{i+1}^{(n)}(r_i^{(n)})\bigr)=0
\end{equation}
for all $i=1,\dots,L$.

By the uniform convergence, $f_1^{(n_j)}(s_{n_j}) \to f_1(s) = x$ and
$y_n^* = f_{L+1}^{(n_j)}(t_{n_j}) \to f_{L+1}(t)$
where we denote the limit by $y^*=f_{L+1}(t)$.
Since $y_n^* \geq y_n$ for all $n$, the inequality remains true in the limit and $y^* \geq y$.
Applying the assumed epigraph lower bound \eqref{epiconv1} for $D_n$ with the convergent endpoints $(x_n,s_n)$ and $(y_n^*,t_n)$ gives
\begin{equation}
\wt D(x,s;y,t)\le D(x,s;y^*,t)\le\liminf_{n\to\infty}D_n(x_n,s_n;y_n^*,t_n)=L,
\end{equation}
which proves property \eqref{epiconv1} for $\wt D$.
\end{proof}

\begin{proof}[Proof of Theorem~\ref{thm:SJtosBr}]
We work under the coupling constructed in Section 4 of \cite{VV26}, under which the rescaled random walk webs $C_n$ converge to the Brownian web $C$ almost surely and $D^\RW_n$ converges to $D^\Br$ almost surely in the epigraph topology.
Here $D^\RW_n$ is the random walk web distance rescaled in the same way as in \eqref{eq:diffusivescaling}.
A lattice path of cost $k$ consists of at most $k+1$ web-path pieces and $k$ cost-one edges, which results in a misalignment level $\varepsilon_n=1/\sqrt n$ for the rescaled distance $D^\RW_n$.
Thus Theorem~\ref{thm:convtransform} applies.
By Proposition \ref{prop:SJinduced}, $D^\SJ_n$ agrees with the half-line version of $D^\RW_n$ provided that the endpoint is reachable from the starting point on some lattice path.
In every compact subset of $\{s<t\}$, the reachability condition is satisfied for all sufficiently large $n$.
Together with \eqref{defDsBr}, this proves convergence on $\{s<t\}$.

It remains to check the epigraph conditions of Lemma \ref{lemma:epigraphconv} with $f_n=D^\SJ_n$ and $f=\mathcal B$ for $s\ge t$.
If $s=t$ and $y\le x$ and $\mathcal B(x,s;y,t)=0$ then we choose $(x_n,s_n)$ for \eqref{epiconv2} as lattice approximations to $(x,s)$ and we follow only up-left edges until $(y_n,t_n)$ where the spatial coordinate approximates $y$ (with a single up-left edge for $y=x$).
These paths have rescaled duration $O(n^{-1/2})$ and the required zero SJ cost.
The lower bound in \eqref{epiconv1} is automatic by nonnegativity.

For $s=t$ and $y>x$ we have $\mathcal B(x,s;y,t)=\infty$ hence \eqref{epiconv2} imposes no condition.
To verify \eqref{epiconv1}, we suppose that a sequence of convergent endpoints has SJ costs bounded by $k$.
Given the corresponding paths $\pi_n$ with SJ cost at most $k$,
there are random walk web paths with at most $k$ cost-one edges and terminal coordinates at least $y_n$ by the first path construction in the proof of Proposition~\ref{prop:SJinduced}.
Stop them at the first lattice vertex at or to the right of a fixed spatial level $z\in(x,y)$; the stopped spatial coordinates tend to $z$.
The finite-piece compactness argument in the proof of Theorem~\ref{thm:convtransform} applies also when the time intervals collapse: each limiting web-path piece is evaluated at the same time $s$, and $\varepsilon_n\to0$ identify all successive endpoint values.
It follows that $z=x$, a contradiction.
Thus the SJ costs tend to infinity and \eqref{epiconv1} holds.

For $s>t$, the cost of any SJ or Brownian web landscape path is infinite by time orientation.
\end{proof}

\section{Particle system representations}\label{s:ips}

This section is devoted to the proofs of Propositions \ref{prop:dTASEP} and \ref{prop:SeppIPS} and Corollary \ref{cor:dTASEPconv}.
The two propositions provide a natural direct connection between the Sepp\"al\"ainen--Johansson model and four interacting particle systems.
The limit theorems in Corollary \ref{cor:dTASEPconv} follow from these relations and from Theorem \ref{thm:SJtosBr}.

For the proof of Proposition \ref{prop:dTASEP}, we first describe how the sequence $(F(0,0;i,j))_{j=0}^\infty$ changes as we increase $i$.
First, observe that $F(0,0;i,j)$ is non-decreasing in $i$ and non-increasing in $j$ by definition, that is,
\begin{equation}\label{Fmonotonicity}
F(0,0;i+1,j)\ge F(0,0;i,j)\ge F(0,0;i,j+1).
\end{equation}

\begin{lemma}\label{lemma:Fevol}
Given $(i,j)$, we let $k$ be the nonnegative integer given by
\begin{equation}\label{defk}
F(0,0;i,j-k-1)>F(0,0;i,j-k)=F(0,0;i,j-k+1)=\dots=F(0,0;i,j).
\end{equation}
Then we have that
\begin{equation}\label{Fevol}
F(0,0;i+1,j)=F(0,0;i,j)+\xi_{i,j-k}\cdot\xi_{i,j-k+1}\cdot\ldots\cdot\xi_{i,j}
\end{equation}
where the number of $\xi$s that we multiply is $k+1$.
\end{lemma}

\begin{proof}[Proof of Lemma \ref{lemma:Fevol}]
The definition \eqref{defF} of the first passage value can be rewritten as the recursion
\begin{equation}\label{Frecursion}\begin{aligned}
F(0,0;i+1,j)&=\min(F(0,0;i,j)+\xi_{ij},F(0,0;i+1,j-1))\\
&=\min\left(\min_{l\in\{0,\dots,m\}}(F(0,0;i,j-l)+\xi_{i,j-l}),F(0,0;i+1,j-m-1)\right)
\end{aligned}\end{equation}
where the second equality above follows by iterating the first equality $m$ times.
If one applies the recursion \eqref{Frecursion} for $m=k$ with $k$ defined in \eqref{defk}
then all the passage values $F(0,0;i,j-l)$ are equal to $F(0,0;i,j)$ for $l=0,\dots,k$ and $F(0,0;i+1,j-k-1)\ge F(0,0;i,j-k-1)>F(0,0;i,j)$.
Hence the minimum on the right-hand side of \eqref{Frecursion} comes from $\min_{l\in\{0,\dots,k\}}(F(0,0;i,j-l)+\xi_{i,j-l})$
and its value is exactly equal to the right-hand side of \eqref{Fevol}.
\end{proof}

\begin{proof}[Proof of Proposition \ref{prop:dTASEP}]
\begin{enumerate}
\item
We consider \eqref{defxij} as the definition of particle positions and we verify that the evolution rules follow discrete-time TASEP.
The evolution rule of passage times given in Lemma \ref{lemma:Fevol} in terms of the particle locations $x_j(i)$ are as follows.
We assume that for given $i$ and $j$ the $j$th particle at time $i$ is in a block of consecutive occupied sites where there are $k$ particles to the right of it, that is,
$x_{j-l}(i)=x_j(i)+l$ for all $l=1,\dots,k$ and $x_{j-k-1}(i)>x_j(i)+k+1$.
In terms of the passage times, this is equivalent to \eqref{defk}.
Hence the $j$th particle moves exactly if $\xi_{i,j-k}=\dots=\xi_{i,j}=1$, that is,
if all particles to the right of the $j$th one in the same block do jump.
This is discrete-time TASEP with sequential update spatially from right to left, that is, increasingly in $j$.

\item
The particle positions can be described in terms of the passage times as
\begin{equation}\label{yjikrewrite}
\{y_j(i)=k-j\}=\{F^{-1}(0,0;i,j)=k\}=\{F(0,0;i,k)\le j,F(0,0;i,k-1)\ge j+1\}
\end{equation}
by \eqref{defyij}.
We assume that $y_j(i)=k-j$ and $y_{j-1}(i)=k+m-(j-1)$ for some $m\ge1$
which means that the gap between particles $j$ and $j-1$ is $m$.
This in particular implies that $F(0,0;i,k)=\dots=F(0,0;i,k+m-1)=j$.
Similarly to \eqref{yjikrewrite}, one can equivalently rewrite the event that the $j$th particle jumps by at least $l$ as
\begin{equation}\label{jumpl}
\{y_j(i+1)\ge k+l-j\}=\{F(0,0;i+1,k+l-1)\ge j+1\}.
\end{equation}
Next, we have that
\begin{multline}
F(0,0;i+1,k+l-1)\\
=\min\left(\min_{r\in\{0,\dots,l-1\}}(F(0,0;i,k+r)+\xi_{i,k+r}),F(0,0;i+1,k-1)\right)
\end{multline}
using \eqref{Frecursion} where $F(0,0;i+1,k-1)\ge F(0,0;i,k-1)\ge j+1$ follows by monotonicity \eqref{Fmonotonicity}
and by assumption \eqref{yjikrewrite}.
Also by assumption, we have $F(0,0;i,k)=\dots=F(0,0;i,k+m-1)=j$ but $F(0,0;i,k+m)\le j-1$.
Hence, \eqref{jumpl} holds and the $j$th particle jumps by at least $l$ if and only if $l\le m$ and $\xi_{i,k}=\dots=\xi_{i,k+l-1}=1$,
that is, with probability $2^{-l}$ for $l=1,\dots,m$.
By $F(0,0;i,k-1)\ge j+1$, the jump is almost surely non-negative and it cannot be larger than $m$.
This corresponds to a discrete-time TASEP with spatial left-to-right sequential update.

\item
The statement is clear in the initial state for $i=0$.
By induction on $i$, we consider a block of particles $x_j(i),\dots,x_{j+l}(i)$ occupying consecutive positions with a hole on their right.
This corresponds to $F(0,0;i,j-1)>F(0,0;i,j)=F(0,0;i,j+l)=k$ for some integer $k$ which means that $x_j(i)=k-j$.
It also implies that $x_j$ has already made $k$ steps by time $i$ and therefore there are exactly $k$ holes behind $x_j(i)$
and the $(k+1)$st hole is on the right of $x_j(i)$ at position $k-j+1$.
In the other system we have by \eqref{yjikrewrite} that $y_k(i)+k=j$ and $y_{k-1}(i)+k-1\ge j+l+1$ meaning that $y_{k-1}(i)-y_k(i)\ge l+2$,
that is, we have a gap of size at least $l+1$ in front of the particle $y_k$.

If $\xi_{ij}=\dots=\xi_{i,j+l}=1$ then all the particles $x_j,\dots,x_{j+l}$ jump to the right by $1$.
On the other hand, the particle $y_k$ in the other system jumps to the right by at least $l+1$ in the time step $i\to i+1$.
Applying the argument for $l=0,1,\dots$ yields a coupling between the $(k+1)$st hole in the first system and the particle $y_k$ in the second system for all $k=0,1,2,\dots$ simultaneously which completes the proof.
\end{enumerate}
\end{proof}

\begin{proof}[Proof of Corollary \ref{cor:dTASEPconv}]
\begin{enumerate}
\item
For $k\ge0$, let $r_{k,n}(t)$ be the rescaled rightmost point reachable from $(0,0)$ with random walk web distance at most $k$.
The reflected-boundary argument of Proposition 4.2 in \cite{VV26}, applied in the forward orientation, gives
\begin{equation}\label{rknconv}
r_{k,n}\stackrel\d\Longrightarrow\rho_k
\end{equation}
uniformly on compact time intervals, where $\rho_k(t)=\sup\{z:\mathcal B(0,0;z,t)\le k\}$.
Indeed, the discrete boundaries have the successive Skorokhod reflection representation; Donsker's theorem and continuity of the reflection map give convergence by induction on $k$.
The one-time law of $\rho_k(1)$ is that of the largest eigenvalue of a $(k+1)$-dimensional GUE random matrix, hence has no atoms; see \cite{B01,W07}.

Let $(y_n,t_n)$ be a deterministic sequence of lattice points rescaled diffusively according to \eqref{eq:diffusivescaling}, that is, $(n^{1/2}y_n,nt_n)\in\Z^2_\e$.
Assume the convergence $(y_n,t_n)\to(y,1)$.
For sufficiently large $n$, the point $(y_n,t_n)$ is reachable from $(0,0)$ in $Z^2_\e$ and the representation \eqref{eq:SJhalfline} together with the definition of $r_{k,n}$ gives that for each $k\ge0$,
\begin{multline}
P\bigl(D_n^\SJ(0,0;y_n,t_n)\le k\bigr)
=P\bigl(r_{k,n}(t_n)\ge y_n\bigr)\\
\to P\bigl(\rho_k(1)\ge y\bigr)
=P\bigl(\mathcal B(0,0;y,1)\le k\bigr)
\end{multline}
where the convergence follows by \eqref{rknconv} and the last equality by the definition of $\rho_k$.
This proves convergence in distribution along such sequences.

Put $m_n=\lfloor n-\sqrt{2n}\,y\rfloor$.
By \eqref{SJF} and \eqref{defxij},
\begin{equation}
x_{m_n}(n)-x_{m_n}(0)=F(0,0;n,m_n)=D^\SJ(0,0;n-m_n,n+m_n).
\end{equation}
With scaling parameter $2n$, the terminal coordinates satisfy
\begin{equation}
\frac{n-m_n}{\sqrt{2n}}\to y,\qquad
\frac{n+m_n}{2n}\to1.
\end{equation}
The preceding convergence proves \eqref{dTASEPconv}.
For $k=0$, $\rho_0(1)$ is standard normal, which gives \eqref{dTASEP0}.

\item
An index $j$ contributes to $N(i,k)$ if and only if $y_j(i)\ge k-j$.
By \eqref{yjikrewrite}, this condition can be written as $F(0,0;i,k-1)\ge j+1$.
Hence we have that $N(n,k)=F(0,0;n,k-1)$.
Taking $k=m_n$ as above gives
\begin{equation}
N(n,m_n)=F(0,0;n,m_n-1)=D^\SJ(0,0;n-m_n+1,n+m_n-1).
\end{equation}
After scaling by $2n$, these terminal coordinates also converge to $(y,1)$, so the first part proves \eqref{dTASEPflux}.
\end{enumerate}
\end{proof}

\begin{proof}[Proof of Proposition \ref{prop:SeppIPS}]
\begin{enumerate}
\item
For the evolution of the system $(w_i(j))$, we let $k=w_i(j-1)$ and $l=w_{i-1}(j)$.
We first observe that by definition \eqref{defwij}
\begin{equation}\label{wwithF}\begin{aligned}
F(0,0;k,j-1)&=i,&F(0,0;k+1,j-1)&=i+1,\\
F(0,0;l,j)&=i-1,&F(0,0;l+1,j)&=i.
\end{aligned}\end{equation}
Since $F$ is non-decreasing in the third variable, by \eqref{defM} which can be written as $M=\max(k,l+1)$, we have that
\begin{equation}\label{FMj}
F(0,0;M,j)=\max(F(0,0;k,j),F(0,0;l+1,j))=i
\end{equation}
because $F(0,0;k,j)\le F(0,0;k,j-1)=i$ by the non-increasing property in the fourth variable and from the last equality in \eqref{wwithF}.
On the other hand, by the first equality in \eqref{Frecursion} for $M+1$, we can write
\begin{equation}\label{FM+1j}
F(0,0;M+1,j)=\min(F(0,0;M,j)+\xi_{M,j},F(0,0;M+1,j-1))=i+\xi_{M,j}
\end{equation}
by \eqref{FMj} and using the fact that $F(0,0;M+1,j-1)\ge F(0,0;k,j-1)=i+1$.
This shows that $w_i(j)\ge M+1$ if and only if $\xi_{M,j}=0$ which is \eqref{wevolution} with $m=1$.
The general case follows by writing \eqref{FM+1j} with $M+1$ replaced by $M+m$ and by induction.

\item
The definition \eqref{defzij} can be rewritten in terms of the first passage values of the Sepp\"al\"ainen--Johansson model as
\begin{equation}
z_i(j)=F(0,0;i+m_i(j),j)+i=m_i(j)+i
\end{equation}
for $i=0,1,2,\dots$ and $j=0,1,2,\dots$ where
\begin{equation}\label{defmij}
m_i(j)=\min\{m\in\Z:F(0,0;i+m,j)\le m\}.
\end{equation}

For the proof of the update rule \eqref{zijevo}, we assume that
\begin{equation}\label{mijassumptions}
m_{i-1}(j)=k,\qquad m_i(j-1)=l
\end{equation}
which corresponds to $z_{i-1}(j)=k+i-1$ and $z_i(j-1)=l+i$.
Since $i\mapsto F(0,0;i,j)$ is nondecreasing, there is equality for the minimal $m$ on the right-hand side of \eqref{defmij}.
In terms of the passage values, \eqref{mijassumptions} can be written as
\begin{equation}\label{Fassumptions}\begin{aligned}
F(0,0;i+k-2,j)&=F(0,0;i+k-1,j)=k,\\
F(0,0;i+l-1,j-1)&=F(0,0;i+l,j-1)=l.
\end{aligned}\end{equation}
Also we assume that there are $r$ horizontal edges with weight $1$ after $z_{i-1}(j)$ before the first free horizontal edge,
that is,
\begin{equation}\label{xiassumptions}
\xi_{i+k-1,j}=\dots=\xi_{i+k+r-2,j}=1,\qquad\xi_{i+k+r-1,j}=0.
\end{equation}
We prove that
\begin{equation}\label{ztildemain}
z_i(j)=\min(z_{i-1}(j)+r+1,z_i(j-1))=\min(i+k+r,i+l)
\end{equation}
in two steps depending on the relation of $k+r$ and $l$.

First let $k+r\le l$ and we prove that $z_i(j)=i+k+r$ by showing that $m_i(j)=k+r$.
For this we use the general passage value identity for the Sepp\"al\"ainen--Johansson model similar to \eqref{Frecursion}
\begin{multline}\label{Frecursion2}
F(0,0;i+k+r-1,j)=\min\Big(F(0,0;i+k-1,j)+\xi_{i+k-1,j}+\dots+\xi_{i+k+r-2,j},\\
\min_{s\in\{1,\dots,r\}}(F(0,0;i+k+s-1,j-1)+\xi_{i+k+s-1,j}+\dots+\xi_{i+k+r-2,j})\Big).
\end{multline}
Since $F(0,0;i+k-1,j)=k$ by \eqref{Fassumptions} and $\xi_{i+k-1,j}=\dots=\xi_{i+k+r-2,j}=1$ by \eqref{xiassumptions},
the first expression in the minimum on the right-hand side of \eqref{Frecursion2} is equal to $k+r$.

On the other hand, for any integers $u\le v$, the inequality
\begin{equation}\label{Fincreasbound}
F(0,0;u,j)\ge F(0,0;v,j)-(v-u)
\end{equation}
follows from the fact that $i\mapsto F(0,0;i,j)$ can increase by at most $1$ in each step.
Since $F(0,0;i+l-1,j-1)=l$ by \eqref{Fassumptions},
applying \eqref{Fincreasbound} with $u=i+k+s-1$ and $v=i+l-1$ means that $F(0,0;i+k+s-1,j-1)\ge k+s$.
We also have that $\xi_{i+k+s-1,j}+\dots+\xi_{i+k+r-2,j}=r-s$, hence each term in the inner minimum on the right-hand side of \eqref{Frecursion2}
is lower bounded by $k+r$ proving $F(0,0;i+k+r-1,j)=k+r$.
Since $\xi_{i+k+r-1,j}=0$, we also have $F(0,0;i+k+r,j)=k+r$ and $m_i(j)=k+r$ and $z_i(j)=i+k+r$.

Next we turn to the case $l<k+r$ and we prove $m_i(j)=l$.
If $l\le k$ then by the first equality in \eqref{Frecursion}, we can write
\begin{equation}\label{Fl<k}
F(0,0;i+l-1,j)=\min(F(0,0;i+l-2,j)+\xi_{i+l-2,j},F(0,0;i+l-1,j-1)).
\end{equation}
By \eqref{Fassumptions}, we have that $F(0,0;i+l-1,j-1)=l$.
On the other hand, the minimality of $k$ in the definition of $m_{i-1}(j)=k$ implies that $F(0,0;i+m-1,j)\ge m+1$ for all $m=k-1,k-2,\dots$.
Since $l\le k$, it is also true for $m=l-1$ which gives that $F(0,0;i+l-2,j)\ge l$.
These together show that the minimum in \eqref{Fl<k} and hence the value of $F(0,0;i+l-1,j)$ is equal to $l$.
By writing
\begin{equation}\label{Fi+l-1toFi+l}\begin{aligned}
F(0,0;i+l,j)&=\min(F(0,0;i+l-1,j)+\xi_{i+l-1,j},F(0,0;i+l,j-1))\\
&=\min(l+\xi_{i+l-1,j},l)\\
&=l
\end{aligned}\end{equation}
using \eqref{Fassumptions}.
This proves $m_i(j)=l$ as required.

If $k<l<k+r$, then we can write $l=k+s$ for some $s\in\{1,\dots,r-1\}$.
The following analog of \eqref{Frecursion2} holds
\begin{multline}\label{Frecursion3}
F(0,0;i+k+s-1,j)=\min\Big(F(0,0;i+k-1,j)+\xi_{i+k-1,j}+\dots+\xi_{i+k+s-2,j},\\
\min_{t\in\{1,\dots,s\}}(F(0,0;i+k+t-1,j-1)+\xi_{i+k+t-1,j}+\dots+\xi_{i+k+s-2,j})\Big).
\end{multline}
For the first expression on the right-hand side, \eqref{Fassumptions} yields $F(0,0;i+k-1,j)=k$ and
$\xi_{i+k-1,j}+\dots+\xi_{i+k+s-2,j}=s$ follows by \eqref{xiassumptions}.
Using \eqref{Fincreasbound} with $v=i+l-1$ and $u=i+k+t-1$, we get the lower bound
$F(0,0;i+k+t-1,j-1)\ge k+t$ since $F(0,0;i+l-1,j-1)=l$ by \eqref{Fassumptions}
and $\xi_{i+k+t-1,j}+\dots+\xi_{i+k+s-2,j}=s-t$ again by \eqref{xiassumptions}.
All this implies that
\begin{equation}
F(0,0;i+l-1,j)=F(0,0;i+k+s-1,j)=k+s=l.
\end{equation}
Then the same argument as in \eqref{Fi+l-1toFi+l} implies $F(0,0;i+l,j)=l$ and $m_i(j)=l$ meaning $z_i(j)=i+l$ which proves \eqref{ztildemain}.

\item
Fix $j$ and write $f(z)=F(0,0;z,j)$ for integers $z\ge0$.
Since $f(0)=0$ and $f(z)-f(z-1)\in\{0,1\}$, the function $z-f(z)$ also has increments in $\{0,1\}$.
By \eqref{defzij}, a site $z\ge1$ is occupied precisely when $z-f(z)>(z-1)-f(z-1)$, or equivalently when $f(z)=f(z-1)$.
Thus $z$ is a hole precisely when $f(z)=f(z-1)+1$.
Putting $i=f(z-1)$, monotonicity of $f$ and \eqref{defwij} give $w_i(j)=z-1$.
Conversely, $w_i(j)=z-1$ implies $f(z-1)=i$ and $f(z)=i+1$, so $z$ is a hole.
All sites $z\le0$ are occupied by convention. This proves the claim.
\end{enumerate}
\end{proof}

\section{The dual symmetry of the Brownian web landscape}
\label{s:duality}

We prove Proposition \ref{prop:duality} in this section about the dual symmetry of the Brownian web landscape.
The natural symmetry of the Sepp\"al\"ainen--Johansson model implies that
\begin{equation}
D^\SJ(i,j;k,l)=\wh D^\SJ(k,l+1;i,j+1)
\end{equation}
for all $(i,j),(k,l)\in\Z^2_\e$ where $\wh D^\SJ$ is the Sepp\"al\"ainen--Johansson distance defined on the dual odd integer lattice.
This symmetry, together with Theorem \ref{thm:SJtosBr} already implies Proposition \ref{prop:duality}.
We give a simple direct proof of the dual symmetry below.

\begin{proof}[Proof of Proposition~\ref{prop:duality}]
We assume that $\mathcal B(x,s;y,t)=k<\infty$.
This means by definition that there exists a path $\pi$ with finitely many jumps
at $(x_1,t_1),\dots,(x_k,t_k)$ between $(x,s)$ and $[y,\infty)\times\{t\}$
where the points of jump are type $(1,2)$ point in the Brownian web $B$ as well as in its dual $\wh B$.
For the first equality in \eqref{duality},
we show the existence of a dual path $\wh\pi$ from $(y,t)$ to $(-\infty,x]\times\{s\}$ with $k$ jumps.

The path $\pi$ must be optimal in the following sense.
For all $i=1,\dots,k$, $\pi$ follows the right outgoing paths starting at the $(1,2)$ point $(x_i,t_i)$
which is to the right of the path $B_{(x_i,t_i)}$ at time $t_{i+1}$,
that is, $\pi(r)=B_{(x_i,t_i)+}(r)$ for $r\in[t_i,t_{i+1}]$ and $B_{(x_i,t_i)}(t_{i+1})<x_{i+1}$ with $(x_{k+1},t_{k+1})=(y,t)$.
If this property did not hold at some point of jump $(x_i,t_i)$ then $[y,\infty)\times\{t\}$ could be reached
with fewer jumps by following the continuation of the incoming Brownian web path at $(x_i,t_i)$ without a jump.

Now we construct the backward path $\wh\pi$.
It starts at $(y,t)$ and it follows first the backward Brownian web path $\wh B_{(y,t)-}$.
By the properties of $\pi$ described above, $B_{(x_k,t_k)}(t)<y\le B_{(x_k,t_k)+}(t)$ holds, hence $\wh B_{(y,t)-}$ must be between the two forward outgoing paths $B_{(x_k,t_k)\pm}$ on $[t_k,t]$ and in particular it goes through $(x_k,t_k)$.
The outgoing backward path $\wh B_{(x_k,t_k)-}$ is to the left of the incoming forward path at $(x_k,t_k)$.
We let $\wh\pi$ follow $\wh B_{(x_k,t_k)-}$ from $(x_k,t_k)$.
By the properties of $\pi$ above, $\wh B_{(x_k,t_k)-}$ starts between $B_{(x_{k-1},t_{k-1})\pm}$
starting at $(x_{k-1},t_{k-1})$ hence $\wh B_{(x_k,t_k)-}$ passes through $(x_{k-1},t_{k-1})$.
By following this procedure inductively, we can define the backward path $\wh\pi$ starting at $(y,t)$ and hitting $(-\infty,x]\times\{s\}$
with jumps exactly at $(x_k,t_k),\dots,(x_1,t_1)$ and following backward Brownian web paths in between.

This construction shows the inequality
\begin{equation}
\mathcal B(x,s;y,t)\ge\wh D^\Br(y,t;(-\infty,x],s)
\end{equation}
jointly for all $(x,s,y,t)\in\R^4$.
The other inequality follows from exactly the same argument.
The equality in distribution
\begin{equation}\label{eq_distr1}
\wh D^\Br(y,t;(-\infty,x],s)\stackrel{\d}=D^\Br(y,-t;(-\infty,x],-s)
\end{equation}
holds by the duality of the forward and backward lines proved in Theorem 2.3 of~\cite{TW98}.
The obvious symmetry of the Brownian web with respect to the reflection about the vertical axis yields
the equality in distribution of the right-hand side of \eqref{eq_distr1} and that of \eqref{duality} for all $y,t,x,s\in\R$ simultaneously.
\end{proof}

\section{Large-scale behavior}\label{s:large-scale}

In this section we prove Theorem \ref{t:large-dev} and
Corollary \ref{cor:LLN}.

The following theorem is from \cite{LedouxRider2010}. The formulation here follows \cite{dauvergne2021bulk}, Theorem 3.1.
\begin{theorem}\label{t:LR}
There exist positive constants $c,d>0$ such that the following holds. Let $k\ge 1$ be an integer, let $t\in[1/2,2]$, and let $\rho_{k-1}$ be the top line of a $k$-level Dyson's Brownian motion started with all eigenvalues at zero and with variance-one Brownian normalization. Then for all $m\in\bigl[0,5k^{2/3}\bigr]$ we have
\begin{equation}\label{e:LR1}
{P}\bigl(\rho_{k-1}(t)-2\sqrt{kt}\geq mk^{-1/6}\bigr)\leq c\exp(-d m^{3/2}),
\end{equation}
and
\begin{equation}\label{e:LR2}    
{P}\bigl(\rho_{k-1}(t)-2\sqrt{kt}\leq -mk^{-1/6}\bigr)\leq c\exp(-d m^3).
\end{equation}
Also, for all $m\geq 5k^{2/3}$  we have
\begin{equation}\label{e:LR3}
{P}\bigl(|\rho_{k-1}(t)-2\sqrt{kt}|\geq mk^{-1/6}\bigr)\leq c\exp(-dk^{-1/3}m^2).
\end{equation}
Note that the upper bound on $m\le 5k^{2/3}$ for \eqref{e:LR1} and 
\eqref{e:LR2} can be replaced by $c'k^{2/3}$ for arbitrary $c'$, simply by using the $5k^{2/3}$ bound when $m\in[5k^{2/3}, c'k^{2/3}]$; this only changes the value of $d$ in the bounds.
\end{theorem}
Let $u=(x_1,s_1;x_2,s_2)$ with $s_1<s_2$ and $\mathcal D(u)\ge 20$. We will cover the domain of $\mathcal B$ where $\mathcal D\ge20$ with rectangular tiles of the form 
$$
T_u=T_u^1\times\dots \times T_u^4, \quad T_u^1=\left\{x:\frac{|x_1-x|}{x_2-x_1}\le \frac{1}{\mathcal D(u)^{2/3}}\right\}, \quad 
T_u^2=\left\{s:\frac{|s_1-s|}{s_2-s_1}\le\frac{1}{ \mathcal D(u)^{2/3}}\right\}.
$$
and $T_u^3,T_u^4$ are defined similarly. 

The next lemma shows that the main term $\D$ is essentially constant over $T_u$.

\begin{lemma}\label{l:D-constant} There exists an absolute constant $c_1$ such that, for all $u$ with $\D(u)\ge 20$, we have 
$$\sup_{v,v'\in T_u}|\D(v)-\D(v')|\le c_1\D(u)^{1/3}.$$
and as a consequence
$$\sup_{v,v'\in T_u}|\D(v)/\D(v')|\le c_1.$$
\end{lemma}
\begin{proof}
Write $u=(x_1,s_1;x_2,s_2)$ and set
$\delta=\D(u)^{-2/3}<1/4$.
For $v=(y_1,t_1;y_2,t_2)\in T_u$, the definition of
the tile gives
$$
1-2\delta
\le \frac{y_2-y_1}{x_2-x_1},
     \frac{t_2-t_1}{s_2-s_1}
\le 1+2\delta.
$$
Therefore 
$$
\frac{(1-2\delta)^2}{1+2\delta}
\le \frac{\D(v)}{\D(u)}
\le \frac{(1+2\delta)^2}{1-2\delta}.
$$
The difference between the upper and lower bounds is
at most $c\delta$, and their ratio is bounded by an
absolute constant. The claims follow since
$\delta\D(u)=\D(u)^{1/3}$.
\end{proof}
\begin{proposition}\label{p:tile}
There exist $c>0$ and $k_0\ge20$ such that, for every $u$ with
$\D(u)\ge k_0$ and every $a\ge0$,
$$
\begin{aligned}
P\left(\sup_{v\in T_u}(\mathcal B-\D)(v)
       \ge a\D(u)^{1/3}\right)
&\le e^{1-ca^2},\\
P\left(\inf_{v\in T_u}(\mathcal B-\D)(v)
       \le -a\D(u)^{1/3}\right)
&\le e^{1-ca^{3/2}}.
\end{aligned}
$$
\end{proposition}
\begin{proof}
{\bf Upper tail bound. } By shift and scale invariance, we may assume that $u=(-x,-1;x,1)$. Then $\mathcal D(u)=\mathbf{1}_{x>0}x^2/2\ge k_0$ so we have $x\ge x_*=\sqrt{2k_0}>0$.

For $(p,q)\in T_u$ with $o=(0,0)$
we use metric composition at the deterministic cutting time $0$ (Proposition~\ref{prop:metriccomp})
$$
\mathcal B(p,q)\le \mathcal B(p,o)+\mathcal B(o,q) 
$$
Set $x'=\sup T_u^3$, $s'=\inf T_u^4$, and $u'=(0,0;x',s')$. Note that $2\D(u')=\D(u'')$ with $u''=(-x',-s';x',s')\in T_u$. We will use the trivial bound
$$
-\D(p,q)\le -2\D(u')+\sup_{v\in T_u}(2\D(u')-\D(v))
$$
Now $T_u$ is defined so that the last supremum is bounded by $c\D(u)^{1/3}$, by  Lemma \ref{l:D-constant}. By Proposition \ref{prop:duality} and the monotonicity of $\mathcal B(o,\cdot,s)$, it suffices to show that (for a possibly  different $d$) with
\begin{equation}\label{e:tile-upper-s-segment}
F=
\left\{\frac{(\sup_{T_u^4}\mathcal B(o;x',\cdot))-\D(u')}{\mathcal D(u)^{1/3}}> a\right\}, \qquad PF\le e^{1-da^2}.
\end{equation}
We first prove the bound at thresholds of the form $$a=\frac{k-\D(u')}{\D(u)^{1/3}},\qquad k>\D(u'),\quad k\in\Z.$$ 
By Proposition \ref{prop:rhoboundary}, the curve $\rho_{k}(t)$ is the left boundary of points of distance at least $k+1$ from $0$. By the reflected Brownian representation in Lemma \ref{lemma:evolutionboundaries} and the Brownian RSK correspondence \cite{B01,W07}, its one-time law is that of the top line of a $(k+1)$-level Dyson's Brownian motion. If 
for some $x_0$ we have $$A=\{x_0\le \rho_{k}(s')\},$$
then all points to the left of the Brownian web path $B_{x_0,s'}$ have distance at most  $k$ from $o$. 

In particular, under the event $F\cap A$ we have $G=\{\inf_{s\in T^4_u}B_{x_0,s'}(s)\le x'\}$. We have
$$
PF\le PA^c+P(F\cap A)\le P(\rho_{k}(s')<x_0)+P
G
$$
Note that by our choice of $u'$ we have $x'=2\sqrt{s'\D(u')}$. Let $x_1=2\sqrt{s'k}$ and set $x_0=(x'+x_1)/2$. Then for $a>0$ we have $$0<x'<x_0=\frac{x'+x_1}{2}<x_1$$
so we can write
$$
A^c=\Big\{\rho_{k}(s')-2\sqrt{s'k}<x_0-x_1\Big\}.
$$
Theorem \ref{t:LR} applies with $k+1$ levels. Its centering $2\sqrt{s'(k+1)}$ is at least $x_1$, and $k^{-1/6}\ge(k+1)^{-1/6}$, which only improve the lower-tail estimates below.
By the reflection principle and Brownian scaling, the random variable $$M=\frac{x_0-\inf_{s\in T^4_u}B_{x_0,s'}(s)}{2\D(u)^{-1/3}}$$ has the distribution of the absolute value of a standard Gaussian, so we get
\begin{equation}\label{e:pg}
PG=P\Big(M\ge (x_0-x')\D(u)^{1/3}/2\Big)\le e^{-(x_0-x')^2\D(u)^{2/3}/8}. \end{equation}

We analyze two cases separately. 

{\bf Case 1, $k\le 2\D(u')$}, where
\begin{equation}\label{e:a-bound}
0<a\le \frac{\mathcal D(u')}{\mathcal D(u)^{1/3}}\le c \mathcal D(u)^{2/3},
\end{equation}
we use the chord lower bound for the concave function $\sqrt{y+z}\ge \sqrt{y}+z/(5\sqrt{y})$ for $z\in[0,y]$. Then  $\sqrt{y+z}- \sqrt{y}\ge z/(5\sqrt{y})$ gives 
\begin{equation}\label{e:x1-x0}\begin{aligned}
x_1-x_0=\frac{x_1
-x'}2
&=
\sqrt{s'}\left(\sqrt{k}-\sqrt{\D(u')}\right)
\ge
c\left(\sqrt{k}-\sqrt{\D(u')}\right)
\\&\ge
c\frac{a\D(u)^{1/3}}{5\sqrt{\D(u')}}
\ge c_0a\D(u')^{-1/6}
\ge c_0ak^{-1/6}
\end{aligned}\end{equation}
for some absolute constant $c_0$ (since $s'$ is bounded from below). We use  Theorem \ref{t:LR}, \eqref{e:LR2}. With that notation, we take $m=c_0a$. The condition $m\le ck^{2/3}$ of the theorem follows from $k\ge\D(u')\asymp\D(u)$, and \eqref{e:a-bound}. Recall $x_1=2\sqrt{s'k}$. We get
\begin{align*}
PA^c\le
P(\rho_{k}(s')-2\sqrt{s'k}<-c_0ak^{-1/6})\le
ce^{-d (c_0a)^3}.
\end{align*}

Next, we use $x_0-x'=x_1-x_0\ge ca \mathcal D(u)^{-1/6}$ by \eqref{e:x1-x0}, and since $\mathcal D(u)$ is bounded below, \eqref{e:pg} gives  $PG\le e^{-ca^2}$.
This concludes the proof of the upper bound in Case 1.

{\bf Case 2, $k\ge 2\D(u')$}, here we use the bound that for $z\ge y$ we have $
\sqrt{y+z}- \sqrt{y}\ge\sqrt{z}/3
$
which gives
\begin{align*}
(x_1-x_0)k^{1/6}&=\frac{1}{2}(x_1-x')k^{1/6}=\sqrt{s'}\left(\sqrt{k}-\sqrt{\D(u')}\right)k^{1/6}
\\&\ge
c\left(\sqrt{k}-\sqrt{\D(u')}\right)k^{1/6}
\ge \frac{c}{3}(a\D(u)^{1/3})^{1/2}(a\D(u)^{1/3})^{1/6}
\\&=c\D(u)^{2/9}a^{2/3}\ge c_1 a^{2/3}
\end{align*}
and in this regime we still have $x_0>0$. We use  Theorem \ref{t:LR}, \eqref{e:LR2}. With that notation, we take $m=c_1a^{2/3}$. The condition $m\le ck^{2/3}$ of the theorem follows from $k\ge a\mathcal D(u)^{1/3}\ge ca$. Recall $x_1=2\sqrt{s'k}$. We get
\begin{align*}
PA^c\le
P(\rho_{k}(s')-2\sqrt{s'k}<-c_1a^{2/3}k^{-1/6})\le ce^{-d a^2}
\end{align*}

Unfortunately, we only get an exponential tail bound for $PG$ in this case. To remedy this, 
for $a\ge1$, partition $T_u^4$ into $\lceil a\rceil$
equal intervals $I_j=[t_j,t_{j+1}]$. Since $t_j\ge s'$
and $x_0>0$, Brownian scaling gives
$$P(\rho_k(t_j)<x_0)
\le P(\rho_k(s')<x_0)
\le c e^{-d a^2}.
$$
Each interval has length $4\D(u)^{-2/3}/\lceil a\rceil$.
Hence the reflection principle gives
$$
P\left(\inf_{s\in I_j}B_{x_0,t_j}(s)\le x'\right)
\le
\exp\left(
-\frac{(x_0-x')^2\lceil a\rceil\D(u)^{2/3}}8
\right)
\le \exp\left(-c a\lceil a\rceil\D(u)\right).
$$
A union bound gives
$$
PF
\le \lceil a\rceil
\left(c e^{-d a^2}
      +e^{-c a\lceil a\rceil\D(u)}\right)
\le c' e^{-d'a^2}.
$$
After adjusting constants and using $PF\le1$ for
bounded $a$, this yields
$
PF\le e^{1-ca^2}
$
completing Case 2 of the upper tail bound. 
Integer rounding changes the threshold by at most one, so adjustment of constants gives the stated bound for every $a\ge0$.

\medskip

\noindent{\bf Lower tail bound.}
We first prove a one-point bound.  By shift and
scale invariance, consider $v=(0,0;x,1)$ with $x>0$.
For an integer $1\le k<\D(v)$, set
$$
a=\frac{\D(v)-k}{\D(v)^{1/3}},
\qquad
z=x-2\sqrt{k}
  =2\left(\sqrt{\D(v)}-\sqrt{k}\right).
$$
By the definition of the boundary $\rho$, we have 
$$
P\bigl(\mathcal B(v)<k\bigr)
=
P\bigl(\rho_{k-1}(1)\ge x\bigr).
$$
Combining the two cases  \eqref{e:LR1} and \eqref{e:LR3} gives for all $z\ge0$
$$
P\bigl(\rho_{k-1}(1)-2\sqrt{k}\ge z\bigr)
\le
c\exp\left(-d\max\{k^{1/4}z^{3/2},z^2\}\right).
$$
If $k\ge\D(v)/2$, then
$$
z=\frac{2(\D(v)-k)}{\sqrt{\D(v)}+\sqrt{k}}
\ge a\D(v)^{-1/6},
\qquad
k^{1/4}z^{3/2}\ge c a^{3/2}.
$$
If $k<\D(v)/2$, then
$$
z^2\ge c\D(v)\ge c a^{3/2}.
$$
Thus $P(\mathcal B(v)<k)\le c e^{-d a^{3/2}}$.
Integer rounding and the trivial probability bound for bounded $a$ give, after adjusting constants, for every $a\ge0$,
\begin{equation}\label{e:point-lower}
P\left(\mathcal B(v)-\D(v)\le-a\D(v)^{1/3}\right)
\le c e^{-d a^{3/2}}.
\end{equation}
The event is empty when its threshold for $\mathcal B(v)$ is negative.

For the uniform bound, normalize $u=(0,0;x,1)$ and
let $r_1=(-x,-1)$ and $r_2=(2x,2)$.
By Proposition~\ref{prop:metriccomp} and $\mathcal B\le\mathcal B_\ev\le\mathcal B+1$, since $r_1$ is deterministic,
$$
\mathcal B(r_1,r_2)\le\mathcal B(r_1,p)+\mathcal B(p,q)+\mathcal B(q,r_2)+2.
$$
Together with Lemma \ref{l:D-constant}, this gives the following bound, absorbing $2$ into $c\D(u)^{1/3}$ since $\D(u)\ge k_0$:
$$
\begin{aligned}
\sup_{(p,q)\in T_u}\bigl(\D(p,q)-\mathcal B(p,q)\bigr)
&\le 3\D(u)-\mathcal B(r_1,r_2)\\
&\quad+
\sup_{p\in T_u^1\times T_u^2}
       \bigl(\mathcal B(r_1,p)-\D(u)\bigr)\\
&\quad+
\sup_{q\in T_u^3\times T_u^4}
       \bigl(\mathcal B(q,r_2)-\D(u)\bigr)
+c\D(u)^{1/3}.
\end{aligned}
$$
The two side rectangles are contained in the tiles
centered at $(r_1;(0,0))$ and $((x,1);r_2)$,
respectively. Both centers have Dirichlet distance
$\D(u)$. Hence the upper-tail estimate already proved,
together with Lemma \ref{l:D-constant}, gives Gaussian
upper tails for the two suprema on the scale
$\D(u)^{1/3}$.

Since $\D(r_1,r_2)=3\D(u)$, the first term is controlled
by \eqref{e:point-lower}. A union bound, followed by
adjustment of constants and the trivial probability
bound for bounded $a$, therefore gives
$$
P\left(
\sup_{(p,q)\in T_u}\bigl(\D(p,q)-\mathcal B(p,q)\bigr)
\ge a\D(u)^{1/3}
\right)
\le e^{1-ca^{3/2}}.
$$

\end{proof}

\begin{proof}[Proof of Theorem \ref{t:large-dev}]
Let $k_0$ be as in the tile Proposition \ref{p:tile}. 
We first consider $\D(u)\ge k_0$. For $j\ge1$, set
$$
K_j=\left\{
u=(x,s;y,t):s<t,\ \D(u)\ge k_0,\
2^j\le 2+\|u\|_\infty+(t-s)^{-1}\le 2^{j+1}
\right\}.
$$
On $K_j$, we have $\|u\|_\infty\le 2^{j+1}$ and
$t-s\ge 2^{-j-1}$. The spatial and temporal radii
of $T_u$ therefore satisfy
$$
(y-x)\D(u)^{-2/3}
=
4^{2/3}(t-s)^{2/3}(y-x)^{-1/3}
\ge c2^{-j}
$$
and
$$
(t-s)\D(u)^{-2/3}
=
4^{2/3}(t-s)^{5/3}(y-x)^{-4/3}
\ge c2^{-3j}.
$$
Partition the box containing $K_j$ into rectangles
of side lengths $c_0 2^{-j}$ in the spatial coordinates
and $c_0 2^{-3j}$ in the temporal coordinates.
For sufficiently small $c_0>0$, choosing one point
of $K_j$ in each rectangle meeting $K_j$ gives a
cover by at most $c2^{12j}$ tiles centered in $K_j$.

By Lemma \ref{l:D-constant}, $\D$ is comparable to its
value at the center for each covering tile.
Also, $\ell(u)$ is comparable to $j$ on $K_j$.
The tile bounds and a union bound give
$$
\begin{aligned}
P\left(
\exists u\in K_j:
(\mathcal B-\D)(u)
\ge a(1+\D(u))^{1/3}\ell(u)^{1/2}
\right)
&\le c2^{12j}e^{-d a^2j},\\
P\left(
\exists u\in K_j:
(\D-\mathcal B)(u)
\ge a(1+\D(u))^{1/3}\ell(u)^{2/3}
\right)
&\le c2^{12j}e^{-d a^{3/2}j}.
\end{aligned}
$$
Summing over $j\ge1$ yields bounds
$c e^{-d a^2}$ and $c e^{-d a^{3/2}}$, respectively,
for all sufficiently large $a$, uniformly over
$\D(u)\ge k_0$.

It remains to consider $\D(u)\le k_0$.
For such $u=(x,s;y,t)$, let
$$
v=\left(x,s;\,x+2\sqrt{k_0(t-s)},t\right).
$$
Then $\D(v)=k_0$, and monotonicity in $y$ gives
$$
(\mathcal B-\D)(u)\le k_0+(\mathcal B-\D)(v).
$$
Moreover, $t-s\le 2\|u\|_\infty$ implies
$\ell(v)\le c\ell(u)$.
Thus, on the event that the upper bound with
coefficient $a$ holds simultaneously for $\D\ge k_0$,
we have
$$
\begin{aligned}
(\mathcal B-\D)(u)
&\le k_0+a(1+k_0)^{1/3}\ell(v)^{1/2}\\
&\le c(1+a)(1+\D(u))^{1/3}\ell(u)^{1/2}.
\end{aligned}
$$
Hence the upper-tail bound extends to all
$s<t$ after adjusting constants.

The lower-tail bound extends to all $s<t$
since, when $\D(u)\le k_0$, deterministically, 
$$
\frac{(\D-\mathcal B)(u)}
     {(1+\D(u))^{1/3}\ell(u)^{2/3}}
\le \frac{k_0}{(\log 2)^{2/3}}.
$$

Define
$$
C_1=
\sup_{\substack{u=(x,s;y,t)\in\R^4\\s<t}}
\frac{(\D(u)-\mathcal B(u))_+}
     {(1+\D(u))^{1/3}\ell(u)^{2/3}},
\qquad
C_2=
\sup_{\substack{u=(x,s;y,t)\in\R^4\\s<t}}
\frac{(\mathcal B(u)-\D(u))_+}
     {(1+\D(u))^{1/3}\ell(u)^{1/2}}.
$$
The estimates above  give the required tail bounds
for sufficiently large $a$, up to fixed prefactors.
Decreasing $c$ and using $P\le1$ for bounded $a$
gives the claimed bounds for every $a\ge0$.
\end{proof}

\section{Structure of geodesics}\label{s:geodesics}

We describe the structure of the geodesics for the distance $\mathcal B$.
An admissible path $\pi$ from $(x,s)$ to $(y,t)$ with $s\le t$ is a c\`adl\`ag function $\pi:[s,t]\to\R$
connecting the halflines $(-\infty,x]\times\{s\}$ and $[y,\infty)\times\{t\}$
so that the curve $u\mapsto(\pi(u),u)$ can only jump to the left in the following sense.
We assume that $\pi(s)\le x$ where the strict inequality corresponds to a jump at the starting time $s$
in which case we use the convention $\pi(s-)=x$.
Further, we have that $\pi(u+)=\pi(u)\le\pi(u-)$ for every $u\in(s,t)$
and that $\pi(t)=y\le\pi(t-)$ where strict inequality is a jump at the terminal time $t$.

A finite subdivision $\tau$ of the interval $[s,t]$ is the finite sequence $s=\tau_0<\tau_1<\dots<\tau_n=t$ for some $n$.
The length of the admissible path $\pi$ is given by
\begin{equation}\label{deflength}
|\pi|=\sup_{\tau=(\tau_0,\dots,\tau_n)}\sum_{i=1}^n \mathcal B(\pi(\tau_{i-1}),\tau_{i-1};\pi(\tau_i-),\tau_i)
\end{equation}
where the supremum is for all finite subdivisions $\tau$ of the interval $[s,t]$.
In the definition \eqref{deflength} above, the distances between the two endpoints of the subintervals over the subdivision $\tau$ are added
but the jumps from $(\pi(\tau_i-),\tau_i)$ to $(\pi(\tau_i),\tau_i)$ do not appear in the formula.
Since $\pi$ is an admissible path, we have $\pi(\tau_i-)\ge\pi(\tau_i)$ and $\mathcal B(\pi(\tau_i-),\tau_i;\pi(\tau_i),\tau_i)=0$ for all $i$.
The path $\pi$ from $(x,s)$ to $(y,t)$ is a geodesic if
\begin{equation}
|\pi|=\mathcal B(x,s;y,t).
\end{equation}

Let $(x,s,y,t)\in\R^4$ such that $s<t$.
Next we describe the structure of all geodesics between $(x,s)$ and $(y,t)$.
We first define two sets of paths $\rho_k$ and $\wh\rho_k$   by induction.
We let $\rho_0(u)=B_{(x,s)+}(u)$ to be the (rightmost) forward Brownian web path starting at $(x,s)$ and
$\wh\rho_0(u)=\wh B_{(y,t)-}(u)$ to be the (leftmost) backward Brownian web path starting at $(y,t)$ for all $u\in[s,t]$.
Then given $\rho_k$ and $\wh\rho_k$, we set
\begin{equation}\label{defrho}
\rho_{k+1}(u)=\sup_{r\in[s,u]}B_{(\rho_k(r),r)+}(u),\qquad
\wh\rho_{k+1}(u)=\inf_{r\in[u,t]}\wh B_{(\wh\rho_k(r),r)-}(u)
\end{equation}
for all $u\in[s,t]$, that is, the supremum (infimum) of all forward (backward) Brownian web paths
starting along the trajectory of the previous path $\rho_k$ ($\wh\rho_k$).

By the results of Subsection 3.2 in~\cite{VV26}, the curves $\wh\rho_k$ are the boundaries of the left neighbourhoods
of the interval $[y,\infty)\times\{t\}$ under the Brownian web distance.
This also holds for the Brownian web landscape according to the next statement.
The Brownian web distance has an irregular behaviour in the second space variable but for the Brownian web landscape
$y\mapsto\mathcal B(x,s;y,t)$ is a nondecreasing function.
Hence an analogous description holds for the forward paths $\rho_k$
which are the boundaries of the right neighbourhoods of $(-\infty,x]\times\{s\}$ under the Brownian web landscape.

\begin{proposition}\label{prop:rhoboundary}
Almost surely, for all $(x,s,y,t)\in\R^4$ with $s<t$, we have that
\begin{align}
\rho_k(u)&=\sup\{z\in\R:\mathcal B(x,s;z,u)\le k\},\label{rhoboundary}\\
\wh\rho_k(u)&=\inf\{z\in\R:\mathcal B(z,u;y,t)\le k\}\label{rhohatboundary}
\end{align}
for all $u\in[s,t]$ and for each $k=0,1,\dots$.
\end{proposition}

\begin{proof}[Proof of Proposition~\ref{prop:rhoboundary}]
The description \eqref{rhohatboundary} immediately follows from \eqref{rhoboundary} by applying the first equality in Proposition~\ref{prop:duality}
to the dual Brownian web.
For the proof of \eqref{rhoboundary}, first observe that $z\mapsto \mathcal B(x,s;z,u)$ is a non-decreasing function.

For any $u\in[s,t]$ and $k=0,1,\dots$, by the coalescence of Brownian web paths,
the supremum in the definition of $\rho_k$ in \eqref{defrho} is achieved.
Hence one can choose a sequence of points $(x_1,t_1),\dots,(x_k,t_k)\in\R^2$ with $t_1<\dots<t_k$ such that
the path following $\rho_i$ between $(x_i,t_i)$ and $(x_{i+1},t_{i+1})$ for $i=0,1,\dots,k$
with $(x_0,t_0)=(x,s)$ and $(x_{k+1},t_{k+1})=(\rho_k(u),u)$ has $k$ jumps which implies $\mathcal B(x,s;\rho_k(u),u)\le k$ and that
\begin{equation}\label{rhoboundaryineq}
\rho_k(u)\le\sup\{z\in\R:\mathcal B(x,s;z,u)\le k\}.
\end{equation}

If there were a strict inequality for some $u$ and $k$ in \eqref{rhoboundaryineq},
then we could consider the pair $(u,k)$ where $u\in[s,t]$ and $k$ is minimal.
For this $u$ and $k$, there must be a path with $k$ jumps from $(x,s)$ to $[z,\infty)\times\{u\}$ for some $z>\rho_k(u)$.
This is impossible for $k=0$.
For $k\ge1$, by the minimality of $k$, the last jump along this path cannot be to the right of $\rho_{k-1}$
and there must be a Brownian web path starting along $\rho_{k-1}$ to $(z,u)$ which contradicts the definition \eqref{defrho}.
\end{proof}

By Proposition 3.5 in \cite{VV26}, the distribution of $\rho_{k+1}$ given $\rho_k$ is equal to that of a Brownian motion reflected off $\rho_k$.
Hence Lemma 7.5 in \cite{VV26} can be used to conclude that $\rho_k(t)\to\infty$ as $k\to\infty$, that is, $\mathcal B(x,s;y,t)$ is finite.
Next we construct the two extremal geodesics from $(x,s)$ to $(y,t)$ using pieces of the boundaries $\rho_k$ with switching times given below.

We define the time of separation for $\rho_k$ and $\rho_{k+1}$ as well as for $\wh\rho_k$ and $\wh\rho_{k+1}$ as
\begin{equation}
\kappa_k(v)=\sup\{u\in[s,v]:\rho_{k+1}(u)=\rho_k(u)\},\quad
\wh\kappa_k(v)=\inf\{u\in[v,t]:\wh\rho_{k+1}(u)=\wh\rho_k(u)\}.
\end{equation}
We remark that the definition of $(\rho_k(\kappa_k(v)),\kappa_k(v))$ is related to that of a relevant separation point for the interval $[s,v]$
in the sense discussed in Subsection 2.3 in~\cite{SSS09} for the Brownian net obtained by the marking construction.

Assume that $\mathcal B(x,s;y,t)=k$.
Then we let $\sigma_{k-1}=\kappa_{k-1}(t)$ and $\wh\sigma_{k-1}=\wh\kappa_{k-1}(s)$.
By induction we further define $\sigma_i=\kappa_i(\sigma_{i+1})$ and $\wh\sigma_i=\wh\kappa_i(\wh\sigma_{i+1})$ for $i=k-2,k-3,\dots,0$.
We let the paths $\pi$ and $\wh\pi$ to be given by
\begin{equation}\label{defpi}
\pi(u)=\left\{\begin{array}{cl} \rho_0(u) & \mbox{if $u\in[s,\sigma_0]$,}\\
\rho_i(u) & \mbox{if $u\in[\sigma_{i-1},\sigma_i]$ for $i=1,\dots,k-1$,}\\
\rho_k(u) & \mbox{if $u\in[\sigma_{k-1},t]$} \end{array}\right.
\end{equation}
and
\begin{equation}\label{defpihat}
\wh\pi(u)=\left\{\begin{array}{cl} \wh\rho_k(u) & \mbox{if $u\in[s,\wh\sigma_{k-1}]$,}\\
\wh\rho_i(u) & \mbox{if $u\in[\wh\sigma_i,\wh\sigma_{i-1}]$ for $i=1,\dots,k-1$,}\\
\wh\rho_0(u) & \mbox{if $u\in[\wh\sigma_0,t]$.} \end{array}\right.
\end{equation}

\begin{remark}
We mention that $\wh\sigma_i\ge\sigma_{k-i-1}$ holds for $i=0,1,\dots,k-1$.
Indeed the reversed inequality would imply the existence of $v\in(\wh\sigma_i,\sigma_{k-i-1})$.
Then the path equal to $\pi$ on $[s,v)$ and to $\wh\pi$ on $[v,t]$ with a left jump from $(\pi(v),v)$ to $(\wh\pi(v),v)$ would have length at most $k-1$ which is impossible.

For the same reason, paths following $\pi$ on $[s,u)$ and then $\wh\pi$ on $[u,t]$ with a jump from $(\pi(u),u)$ to $(\wh\pi(u),u)$ for some $u\in(s,t)$ are not necessarily geodesics because their length can exceed $k$.
\end{remark}

\begin{theorem}\label{thm:geodesicchar}
For $(x,s,y,t)\in\R^4$ with $s<t$ and $\mathcal B(x,s;y,t)=k$, we consider the paths $\pi$ and $\wh\pi$ defined above.
\begin{enumerate}
\item The path $\pi$ defined in \eqref{defpi} is the rightmost geodesic between $(x,s)$ and $(y,t)$.
\item The path $\wh\pi$ defined in \eqref{defpihat} is the leftmost geodesic between $(x,s)$ and $(y,t)$.
\item A c\`adl\`ag path $\gamma$ between $\wh\pi$ and $\pi$ is a geodesic between $(x,s)$ and $(y,t)$
if and only if the following conditions hold.
The path $\gamma$ only has jumps down: $\gamma(u)\le x$, $\gamma(u)=\gamma(u+)\le\gamma(u-)$ for all $u\in(s,t)$, and $\gamma(t)=y\le\gamma(t-)$.
Except for $k$ points, the modulus of continuity of $\gamma$ can locally be upper bounded by that of the Brownian web $B$,
that is, for every nonexceptional $u\in[s,t)$, the inequality $\gamma\le B_{(\gamma(u),u)}$ holds on some interval $[u,u+\varepsilon]$.
At the $k$ exceptional points, the inequality $\gamma\le B_{(\gamma(u),u)+}$ holds locally.
The $k$ exceptional points are exactly the points of increase of the geodesic length $u\mapsto|\gamma_{[s,u]}|$.
\end{enumerate}
\end{theorem}

\begin{remark}
For any admissible path $\gamma$ we define $\wh\gamma$ which differs from $\gamma$ only at the points of jumps $v$
where it is given by $\wh\gamma(v)=\gamma(v-)$ instead of $\gamma(v+)$.
Then the flipped backward path $u\mapsto-\wh\gamma(-u)$ is admissible from $(-y,-t)$ to $(-x,-s)$.

By the first equality in Proposition \ref{prop:duality},
the Brownian web landscape $\wh{\mathcal B}$ defined in terms of the dual Brownian web $\wh B$ is almost surely equal to $\mathcal B$.
This implies that the length of an admissible path $\gamma$ defined with $\mathcal B$ by \eqref{deflength}
is equal to the length of $\wh\gamma$ defined with $\wh{\mathcal B}$ and that $\gamma$ is a geodesic if and only if $\wh\gamma$ is.
Hence a similar characterization can be given as in the third part of Theorem \ref{thm:geodesicchar}
for the geodesic property of a path in terms of its backward modulus of continuity compared to that of the backward Brownian web $\wh B$.
The equivalence of the two descriptions in terms of forward and backward modulus of continuity is not a priori clear.
\end{remark}

\begin{proof}[Proof of Theorem \ref{thm:geodesicchar}]
\begin{enumerate}
\item
$\pi$ is clearly a geodesic.
Assume that $\gamma$ is a geodesic from $(x,s)$ to $(y,t)$ and define
\begin{equation}
\sigma=\inf\{u\in[s,t]:\gamma(u)>\rho_{k-1}(u)\}
\end{equation}
to be the time of separation of $\gamma$ and $\rho_{k-1}$ for which one has $\gamma(\sigma)=\rho_{k-1}(\sigma)$.
By Proposition~\ref{prop:rhoboundary}, $\gamma$ cannot have any more jumps after time $\sigma$
hence it must follow the Brownian web path $B_{(\gamma(\sigma),\sigma)+}$ starting at the $(1,2)$ point $(\gamma(\sigma),\sigma)$,
that is, $\gamma(u)\le\rho_k(u)=\pi(u)$ for $u\in[\sigma,t]$.
By induction, it follows that $\gamma\le\pi$ holds on the full interval $[s,t]$.

\item
Although $\wh\pi$ is not a path with finitely many jumps in the forward Brownian web $B$, it is still a geodesic.
Indeed, its length is computed as the supremum of the sum of the distances between subdivision points
for any finite subdivision of $[s,t]$.
By applying the first equality in Proposition~\ref{prop:duality} for each distance in the subdivision,
one recovers that the length of $\wh\pi$ is $k$
which is equal to its length according to the Brownian web landscape defined on the dual web $\wh B$.
The fact that it is a leftmost geodesic follows as in the previous part.

\item
Assume that the conditions for $\gamma$ are satisfied, then we prove that it is a geodesic by showing that its length is $k$.
By definition \eqref{deflength}, we consider an arbitrary finite subdivision $\tau$ of $[s,t]$.
By refining the subdivision if necessary we have $k$ intervals $[\tau_{i-1},\tau_i]$ where the modulus of continuity of $\gamma$
is not bounded by that of the Brownian web.
In these intervals $\gamma(\tau_i)\le B_{(\gamma(\tau_{i-1}),\tau_{i-1})}(\tau_i)$ fails to hold but
$\mathcal B(\gamma(\tau_{i-1}),\tau_{i-1};\gamma(\tau_i),\tau_i)=1$ is still true
because the Brownian web path $B_{(\gamma(\tau_{i-1}),\tau_{i-1})}$ can be followed
and then use the additional right outgoing path at an appropriate $(1,2)$ point and possibly go to the left for free in the end.
In the other intervals the distance of endpoints is $0$ since the Brownian web can be followed and moving to the left in the end is free.
It means that the length of $\gamma$ is $k$ hence it is a geodesic.

For proving the necessity of the conditions, assume now that $\gamma$ is a geodesic.
The length $u\mapsto|\gamma_{[s,u]}|$ must be a non-decreasing function from $[s,t]$ to $\{0,1,\dots,k\}$.
We show that it takes all the values in $\{0,1,\dots,k\}$, that is, it cannot have a multiple jump.
A double jump would mean that $\rho_k$ for three consecutive values of $k$ would go through the same point $(\gamma(u),u)$
for some $u\in(s,t)$.
Next we prove that it is impossible.

Without loss of generality we can assume that $\rho_0,\rho_1,\rho_2$ go through $(\gamma(u),u)$.
The double jump of the length function also means that $\rho_1$ and $\rho_2$ do not meet
in time $(u,u+\varepsilon)$ for some $\varepsilon>0$,
otherwise the geodesic could follow $\rho_1$ with a single jump instead.
The point $(\gamma(u),u)$ has an incoming path in the Brownian web $B$ and two outgoing path.
The left outgoing path is $\rho_0=B_{(\gamma(u),u)}$.
The right outgoing path $B_{(\gamma(u),u)+}$ which is equal to $\rho_2$ locally and it is disjoint from $\rho_0$ on $(u,u+\varepsilon)$.
Hence $(\gamma(u),u)$ is a type $(1,2)$ point.
By the definition of $\rho_1$ in \eqref{defrho}, we have $\rho_1\ge B_{(\gamma(u),u)+}=\rho_2$ locally which is a contradiction.

This means that there are $k$ distinct times where the length function $u\mapsto|\gamma_{[s,u]}|$ grows by $1$.
Let $v$ be a time of increase.
We show that $(\gamma(v),v)$ is a $(1,2)$ point in the Brownian web $B$.
There are three possibilities for the relation of the paths $B_{(\gamma(v),v)}$ and $\gamma$
on some interval $(v,v+\varepsilon)$ locally.
If $\gamma\le B_{(\gamma(v),v)}$ locally then there cannot be a jump at $v$ in the length function
by Proposition~\ref{prop:rhoboundary}.
If $\gamma>B_{(\gamma(v),v)}$ locally then $\gamma$ must follow the right outgoing path $B_{(\gamma(v),v)+}$
starting at $(\gamma(v),v)$ which means that this point has type $(1,2)$.
If there are times $v'<v''$ arbitrarily close to $v$ such that $\gamma(v')>B_{(\gamma(v),v)}(v')$
and $\gamma(v'')\le B_{(\gamma(v),v)}(v'')$ then the geodesic $\gamma$ could be redefined to follow $B_{(\gamma(v),v)}$
on the time interval $[v,v'']$ which would result in a smaller weight and it would lead to a contradiction.

As a conclusion, all the $k$ jumps of the length function are $(1,2)$ points of $B$
where $\gamma$ follows locally the right outgoing path meaning that the modulus of continuity condition is satisfied
at the $k$ exceptional points.
The modulus of continuity bound at all other points and the condition on the direction of jumps follows from the description
given in Proposition~\ref{prop:rhoboundary} for $k=0$.
\end{enumerate}
\end{proof}

\begin{proof}[Proof of Corollary \ref{cor:infgeod}]
Let $\pi$ be the rightmost geodesic from $(x,s)$ to $(y,t)$ in Theorem \ref{thm:geodesicchar}.
By construction, $\pi$ follows the Brownian web path $B_{(x,s)}$ on the time interval $[s,s+\varepsilon]$ for some $\varepsilon>0$.
Then we show that there are uncountably many further geodesics.

By possibly decreasing $\varepsilon$, we can assume that the point $(z,s+\varepsilon)$ with $z=\pi(s+\varepsilon)$
has a single incoming path in the forward Brownian web,
that is, it is a type $(1,1)$ point in the forward web and a type $(0,2)$ point in the backward web.
We consider the leftmost backward Brownian web path $\wh B_{(z,s+\varepsilon)-}$ starting at $(z,s+\varepsilon)$
which cannot coincide with $\pi$ on $[s,s+\varepsilon]$.
Hence there is $s'\in(s,s+\varepsilon)$ such that $\wh B_{(z,s+\varepsilon)-}(s')<\pi(s')$.

For any $w\in(\wh B_{(z,s+\varepsilon)-}(s'),\pi(s'))$, we construct a geodesic $\pi_w$ which agrees with $\pi$ on $[s,s']$,
jumps to $(w,s')$ at time $s'$ and then follows the rightmost forward Brownian web path $B_{(w,s')+}$ until $s+\varepsilon$
which is well-defined for all $w$.
Since $B_{(w,s')+}$ cannot cross $\wh B_{(z,s+\varepsilon)-}$, $\pi_w$ must go through $(z,s+\varepsilon)$.
The length of $\pi_w$ is $0$ on $[s,s+\varepsilon]$ by construction hence it is a geodesic.
The paths $\pi_w$ for $w\in(\wh B_{(z,s+\varepsilon)-}(s'),\pi(s'))$ are all different and their number is more than countable.
The construction works for all $(x,s,y,t)\in\R^4$ with $s<t$ unless the starting point $(x,s)$ is exceptional with multiple outgoing paths
in which case the rightmost geodesic $\pi$ starts following $B_{(x,s)+}$ and the rest of the argument is the same.
\end{proof}

\section{Height evolution}
\label{s:evolution}

The metric composition is a key property of the directed landscape.
The following version of the Brownian web landscape which is used for the evolution of distances satisfies the metric composition.
Let $\mathcal B_\ev(x,s;y,t)$ denote the minimum cost of a continuous path
from $(x,s)$ to $[y,\infty)\times\{t\}$ that follows finitely many Brownian
web paths. Each switch between web paths in $(s,t)$ costs one.
In addition, starting along the additional outgoing path at a type $(1,2)$
point costs one; all other initial choices cost zero.
Thus
\begin{equation}
\mathcal B(x,s;y,t)\le \mathcal B_\ev(x,s;y,t)\le \mathcal B(x,s;y,t)+1,
\end{equation}
with equality in the first inequality unless $(x,s)$ has type $(1,2)$.
At points with no incoming path, any outgoing path may be chosen without
an initial charge.

Let $h_0:\mathbb R\to \mathbb R\cup \{\infty\}$ be a lower semicontinuous function
that is bounded from below on $\R_+$ and satisfies the following quadratic decay bound at $-\infty$:
\begin{equation}\label{h0decay}
\liminf_{x\to-\infty}\frac{h_0(x)}{x^2}\ge0.
\end{equation}
Consider the evolution
\begin{equation}\label{defevolution}
h_t(y)=\inf_{x\in \mathbb R}\left(h_0(x)+\mathcal B_\ev(x,0;y,t)\right).
\end{equation}

\begin{proposition}\label{prop:metriccomp}
\begin{enumerate}
\item
Almost surely, for all times $r<s<t$ and locations $x,y\in\R$, we have that
\begin{equation}\label{metriccomp}
\mathcal B_\ev(x,r;y,t)=\inf_{z\in\R}\left(\mathcal B_\ev(x,r;z,s)+\mathcal B_\ev(z,s;y,t)\right).
\end{equation}
\item
Almost surely, for all times $s<t$ and locations $y\in\R$, we have that
\begin{equation}\label{htcomp}
h_t(y)=\inf_{z\in\R}\left(h_s(z)+\mathcal B_\ev(z,s;y,t)\right).
\end{equation}
\end{enumerate}
\end{proposition}

\begin{proof}
\begin{enumerate}
\item
Fix $z\in\R$ and choose minimizing paths $\pi_1$ and $\pi_2$ for the two
terms on the right-hand side of \eqref{metriccomp}. Write $w=\pi_1(s)\ge z$.
If $w=z$, concatenate the paths. Any switch at time $s$ is already
charged as an initial switch of $\pi_2$.
If $w>z$, continue the last web path of $\pi_1$ until it first meets
$\pi_2$, and then follow $\pi_2$. At a first meeting before $t$, the two
incoming paths are distinct, so the meeting point has type $(2,1)$ and
joining $\pi_2$ costs nothing. If there is no meeting, the continued path
ends to the right of $\pi_2(t)\ge y$.
In either case the resulting path has cost at most the sum on the
right-hand side of \eqref{metriccomp}. Taking the infimum over $z$ proves
one inequality.

Conversely, split a minimizing path for the left-hand side at time $s$
and set $z=\pi(s)$. Assign any switch at $s$ to the second piece.
The initial-charge convention makes the sum of the costs of the two
pieces equal to the original cost, and each piece costs at least the
corresponding distance in \eqref{metriccomp}. This proves the reverse
inequality.

\item
It follows immediately from \eqref{metriccomp} and the definition of the evolution in \eqref{htcomp}.
\end{enumerate}
\end{proof}

For each deterministic $s$, almost surely there are no type $(1,2)$ points
on $\R\times\{s\}$, by Theorem 3.14 in \cite{FINR06}. Hence
$\mathcal B_\ev(x,s;y,t)=\mathcal B(x,s;y,t)$ simultaneously for all $x,y$ and $t>s$.
Applying \eqref{metriccomp} at the deterministic starting times $0$ and $s$
therefore gives, almost surely for each deterministic $0<s<t$, and $y\in\R$ the definitions \eqref{evolution} and \eqref{defevolution} coincide.

The Brownian web restricted to times after $s$ is independent of its
restriction to times before $s$, and has a time-translation invariant law.
Indeed, construct the past and future webs independently and concatenate
each past skeleton path with the future path from its endpoint at $s$.
Finite families have the coalescing Brownian motion law, so taking the
closure gives the Brownian web. Time-translation invariance follows from
the same construction.
Consequently \eqref{defevolution} is a time-homogeneous Markov evolution
for deterministic initial data, or initial data independent of the web.
The following lemma establishes preservation of the state space.

\begin{lemma}\label{lemma:evolutionprop}
Almost surely, for each $t>0$, the following hold.
\begin{enumerate}
\item The function $y\mapsto h_t(y)$ is nondecreasing.
\item If $h_0$ is lower semicontinuous then so is $h_t$.
\item The function $h_t$ satisfies the quadratic decay bound \eqref{h0decay} with $h_0$ replaced by $h_t$, and is bounded from below on $\R_+$.
\end{enumerate}
\end{lemma}

\begin{proof}[Proof of Lemma \ref{lemma:evolutionprop}]
\begin{enumerate}
\item
As explained before the lemma, it holds almost surely that for $s=0$, $\mathcal B_\ev(x,0;y,t)=\mathcal B(x,0;y,t)$ for all $x,y\in\R$ and all $t>0$.
Note that $y\mapsto \mathcal B(x,0;y,t)$ is nondecreasing by the definition \eqref{defDsBr}.
Hence the function $y\mapsto h_t(y)$ is also nondecreasing as the infimum of nondecreasing functions.

\item
The Brownian web landscape $\mathcal B(x,0;y,t)$ is lower semicontinuous in $y$ by the fact that it is defined using the lower semicontinuous version of the Brownian web distance.
Hence $h_t(y)$ is lower semicontinuous as the infimum of lower semicontinuous functions by definition \eqref{defevolution}.

\item
Fix $t>0$ and $0<\delta<1/(32t)$. The assumptions on $h_0$
imply that $h_0(x)\ge-c_\delta-\delta(x\wedge0)^2$ for all $x$.
For $y<0$ and $x\ge2y$, nonnegativity of the distance therefore gives
$$
h_0(x)+\mathcal B(x,0;y,t)\ge-c_\delta-4\delta y^2.
$$
For $x<2y<0$, we have
$$
\D(x,0;y,t)=\frac{(y-x)^2}{4t}\ge\frac{x^2}{16t}.
$$
Theorem \ref{t:large-dev} implies, uniformly over $x<2y$ as $y\to-\infty$,
that $\mathcal B(x,0;y,t)\ge x^2/(32t)$: its lower-tail error is
$O(|x|^{2/3}(\log(2+|x|))^{2/3})$ in this region.
Hence the contribution from $x<2y$ is bounded below by $-c_\delta$.
Taking the infimum over $x$ gives
$h_t(y)\ge-c_\delta-4\delta y^2$ for all sufficiently negative $y$.
Letting $\delta\downarrow0$ proves \eqref{h0decay}. Monotonicity gives $h_t(y)\ge h_t(0)>-\infty$ for $y\ge0$.

\end{enumerate}
\end{proof}

Moreover, the evolution \eqref{defevolution} started from a given $h_0$
agrees with the evolution started from the backward running minimum process of $h_0$, that is from $x\mapsto\inf_{y\ge x}h_0(y)$.
Thus we may restrict the state space of the Markov process defined by \eqref{defevolution} to nondecreasing lower semicontinuous functions.

Since $\mathcal B$ is integer-valued, the case when $h_0$ has values in $\mathbb Z\cup \{\infty\}$ is of special interest.
For the boundary lemma and Fredholm determinant formula below, we restrict our attention to a further special case when $h_0$ is $\Z\cup\{\infty\}$-valued and bounded from below.

Assume that $h_0$ is nondecreasing, lower semicontinuous and $\Z\cup\{\infty\}$-valued,
with $\lim_{x\to-\infty}h_0(x)=0$ and $\lim_{x\to\infty}h_0(x)=\infty$.
For $i\ge1$ let

\begin{equation}\label{defYti}
Y_t(i)=\max\{x:h_t(x)<i\}.
\end{equation}

\begin{lemma}\label{lemma:evolutionboundaries}
We have $Y_t(1)=B_{(Y_0(1),0)}(t)$.
Given the paths $(Y_s(1),\ldots,Y_s(k-1))_{s\ge0}$,
 the conditional law of the path $(Y_t(k))_{t\ge0}$ is that of Brownian motion reflected in the Skorokhod sense from $Y_t(k-1)$ above.
Equivalently, their joint law is given by independent standard Brownian motions $W_k$ through
\begin{equation}
Y_t(1)=Y_0(1)+W_1(t),\qquad
Y_t(k)=Y_0(k)+W_k(t)+\sup_{0\le r\le t}\bigl(Y_r(k-1)-Y_0(k)-W_k(r)\bigr)_+
\end{equation}
for $k\ge2$.
\end{lemma}

\begin{proof}[Proof of Lemma \ref{lemma:evolutionboundaries}]
The characterization of $Y_t(1)$ follows from the fact that $h_t(x)=0$ holds if and only if there exist $x'\le {Y_0(1)}$ and $x''\ge x$ such that $D^\Br(x',0;x'',t)=0$.
This means that $(x'',t)$ must be on the Brownian web path $B_{(x',0)}$ started at $(x',0)$, hence the boundary curve for points with $h_t(x)=0$ is $B_{(Y_0(1),0)}$.

Next we use induction on $k$.
Given $Y_t(1),\dots,Y_t(k-1)$, $h_t(x)=k-1$ can happen at a point $(x,t)$ in two ways.
Either the Brownian web path $B_{(Y_0(k),0)}$ started at $(Y_0(k),0)$ or the rightmost Brownian web path $B_{(Y_s(k-1),s)+}$ for some $s\in[0,t]$ started at some point along $Y(k-1)$ is weakly to the right of $x$ at time $t$.
This means that the $k$th boundary $Y(k)$ is the supremum of the all Brownian web paths started at $(Y_0(k),0)$ and along the $(k-1)$st boundary.
Next we describe the conditional law of the boundary $Y(k)$.

We reveal $Y(1),\ldots,Y(k-1)$ and put $f=Y(k-1)$.
The one-sided exploration in the proof of Proposition 3.5 in \cite{VV26}, with time orientation reversed, shows that the web strictly to the right of $f$ remains unexplored.
The same exploration applies to ordered initial positions: the new path starts at $Y_0(k)\ge f(0)$ and follows the Brownian web until it first meets $f$.

For $\eps>0$, we approximate the new boundary $Y(k)$ by following the Brownian web path from $(Y_0(k),0)$ until it meets $f$, restarting at height $f+\eps$, and following the Brownian web path.
Repeating these steps defines an approximation $Z^\varepsilon(t)$ for $Y_t(k)$.
Conditionally on the revealed boundaries, each new segment uses fresh Brownian increments from the unexplored side.
Concatenating these increments gives a standard Brownian motion $W_k$ whose conditional law does not depend on $Y(1),\ldots,Y(k-1)$, hence $W_k$ is independent of them.
The driving Brownian motion $W_k$ can be defined by writing the approximating boundary as
\begin{equation}\label{Zeps}
Z^\eps(t)=Y_0(k)+W_k(t)+L^\eps(t),
\end{equation}
where $L^\eps$ increases by $\eps$ whenever $Z^\eps$ meets $f$ (including at time zero if $Y_0(k)=f(0)$).

By definition of the driving Brownian motion $W_k$,
\begin{equation}
L(t)=\sup_{0\le r\le t}\bigl(f(r)-Y_0(k)-W_k(r)\bigr)_+
\end{equation}
satisfies $0\le L^\eps(t)-L(t)\le\eps$.
Thus the approximations in \eqref{Zeps} converge in conditional law as $\varepsilon\to0$, uniformly on compact intervals, to $Y_0(k)+W_k+L$.
The web approximation argument of Proposition 7.3 in \cite{VV26} identifies their limit with the boundary $Y(k)$.
This proves the claimed conditional law and, by induction, the representation with independent driving Brownian motions.
\end{proof}

\begin{theorem}\label{thm:RBM-formula}
Consider the process $h_t$ with nondecreasing lower semicontinuous integer-valued initial condition $h_0$ satisfying $\lim_{x\to -\infty} h_0(x)=0$ and $\lim_{x\to\infty}h_0(x)=\infty$.
For any values $1 \leq n_1 < n_2 < \dotsc <n_m$, any locations $a_1,\dotsc,a_m\in\mathbb R$ and any $t > 0$, we have
\begin{equation}\label{RBM-formula}
P\left(h_t(a_j) \ge n_j, j=1,\dotsc,m \right) = \det\left(\mathbf 1-\bar\chi_{-a} K_t\bar\chi_{-a}\right)_{L^2(\{n_1,\dotsc,n_m\}\times\mathbb R)}
\end{equation}
where $\det$ is the Fredholm determinant with
\begin{multline}\label{KtRBM-Hermite}
K_t(n_i,z_i;n_j,z_j)=-\partial^{-(n_j-n_i)}(z_i,z_j)1_{n_i<n_j}\\
+\sum_{\ell=0}^{n_j-1}\iint_{\mathbb R^2}\d\eta\,e^{b-\eta}P_{\eta}(N(\sigma)=\ell,\,\sigma\in\d b)\,\psi_{n_i}(t,-\eta-z_i)\bar\psi_{n_j-\ell-1}(t,-b-z_j)
\end{multline}
where $P_\eta$ is the law of a rate-one, right-continuous Poisson counting process $N(x)$ started at time $\eta$ with $N(\eta)=0$.
Further, we have $\sigma=\inf\{x\ge\eta:N(x)\ge h_0(x)\}$, with
$\inf\varnothing=\infty$; the hitting measure in \eqref{KtRBM-Hermite}
is restricted to $\sigma<\infty$ and includes the atom at $\sigma=\eta$.
For $n\ge1$, define $\partial^{-n}(z,z')=(z-z')^{n-1}\mathbf1_{z>z'}/(n-1)!$, and let
$H_n(x)=(-1)^n e^{x^2/2}\frac{\d^n}{\d x^n}e^{-x^2/2}$ be the probabilists' Hermite polynomial.
We set
\begin{equation}\label{eq:uppsi}
\psi_n(t,x)=t^{-n/2} \tfrac{1}{\sqrt{2\pi t}}e^{-\frac{x^2}{2t}} H_{n}(\tfrac{x}{\sqrt{t}}),\qquad\bar\psi_n(t,x)=\tfrac{1}{n!}t^{n/2}H_{n}(\tfrac{x}{\sqrt{t}})
\end{equation}
and for $a=(a_1,\dots,a_m)\in\R^m$ and $n_1<\dots<n_m$, we have $\bar\chi_a(n_j,x)=\mathbf1_{x\le a_j}$.
\end{theorem}

\begin{proof}[Proof of Theorem \ref{thm:RBM-formula}]
Let $X_t(i)=-Y_t(i)$ for all $t\ge0$ and integer $i\ge1$.
Then $X_t(i)$ are one-sided reflected Brownian motions with position decreasing with the index as in \cite{NQR20}.
The assumptions on $h_0$ imply that every $Y_0(i)$ is finite and that $X_0(i)=-Y_0(i)$ is nonincreasing in $i$.
Coincident initial positions are allowed in the Skorokhod construction of \cite[Section 1]{NQR20}.
Moreover, the first $n_m$ particles evolve independently of subsequent particles, so this application of Theorem 2.1 requires no additional growth assumption on $X_0(i)$ as $i\to\infty$.

By Theorem 2.1 of \cite{NQR20}, for any $1 \leq n_1 < n_2 < \dotsc <n_m$, $a_1,\dotsc,a_m\in\mathbb R$ and any $t > 0$, we have
\begin{equation}\label{eq:extKernelProbRBM}
P\left(X_t(n_j) \geq a_j, j=1,\dotsc,m \right) = \det\left(\mathbf 1-\bar\chi_a\wt K_t\bar\chi_a\right)_{L^2(\{n_1,\dotsc,n_m\}\times\mathbb R)}
\end{equation}
where
\begin{multline}\label{KtRBMtilde}
\wt K_t(n_i,z_i;n_j,z_j)=-\partial^{-(n_j-n_i)}(z_i,z_j)1_{n_i<n_j}\\
+\sum_{\ell=0}^{n_j-1}\iint_{\mathbb R^2}\d\eta\,e^{\eta-b}P_{B_0=\eta}(\tau=\ell,\,B_\tau\in\d b)\,\psi_{n_i}(t,\eta-z_i)\bar\psi_{n_j-\ell-1}(t,b-z_j)
\end{multline}
where $(B_k)_{k\geq0}$ is a discrete time random walk taking Exp$(1)$ steps to the left, $\tau=\inf\left\{k\geq0\!:B_k\geq X_0(k+1)\right\}$
is the hitting time of the epigraph of the curve $(X_0(k+1))_{k\geq0}$ by $B_k$.

By lower semicontinuity and \eqref{defYti},
\begin{equation}\label{httoXt}\begin{aligned}
\{h_t(a_j)\ge n_j,\ j=1,\dots,m\}
&=\{Y_t(n_j)<a_j,\ j=1,\dots,m\}
\\&=\{X_t(n_j)>-a_j,\ j=1,\dots,m\}.
\end{aligned}\end{equation}
Each $Y_t(i)$ has an atomless distribution for $t>0$. 

To see this, subtract
the average of the first $n_m$ driving Brownian motions from each individual one. 
The average is independent of the centered versions.
The Skorokhod construction commutes with adding the same continuous path to all driving paths.
This follows by induction from a deterministic fact: $f_1+g$ Skorokhod reflected from $f_2+g$ equals  $g+\hat f$, where $\hat f$ is  $f_1$ Skorokhod reflected from $f_2$.

Thus each $Y_t(i)$, $i\le n_m$, is the sum of an independent Gaussian
of variance $t/n_m$ and a random variable determined by the centered
driving Brownian motions. Thus $P(\{Y_i\le a_i\}\setminus\{Y_i<a_i\})=0$ and we may therefore replace $>$ by $\ge$ on the right-hand side of \eqref{httoXt}
and apply \eqref{eq:extKernelProbRBM} with thresholds $-a_j$.

The process
\begin{equation}
N(x)=\max\{k\ge0:-B_k\le x\}
\end{equation}
for $x\ge-\eta$ is a rate-one Poisson counting process with $N(-\eta)=0$.
By \eqref{defYti},
\begin{equation}
\tau=\inf\{k\ge0:-B_k\le Y_0(k+1)\}
=\inf\{k\ge0:h_0(-B_k)\le k\}.
\end{equation}
Since $h_0$ is nondecreasing and $N$ is constant between successive
arrival times, the first occurrence of $N(x)\ge h_0(x)$ is either at
$x=-\eta$ or at an arrival time. Consequently, on $\{\tau<\infty\}$,
\begin{equation}
\sigma=\inf\{x\ge-\eta:N(x)\ge h_0(x)\}=-B_\tau,
\qquad N(\sigma)=\tau,
\end{equation}
and $\sigma=\infty$ if $\tau=\infty$.
Changing variables $\eta\to-\eta$ and $b\to-b$ in
\eqref{KtRBMtilde} now gives \eqref{KtRBM-Hermite}: the factor
$e^{\eta-b}$ becomes $e^{b-\eta}$, the Hermite arguments become
$-\eta-z_i$ and $-b-z_j$, and the index $\ell=\tau=N(\sigma)$ is unchanged.
\end{proof}

\subsection{Increment-stationary laws}
\begin{proof}[{\color{DarkGreen}Proof of Theorem \ref{prop:poissonstationarity}}]
For $k\in\Z$, write
$$
Y_0(k)=\sup\{x:h_0(x)<k\}.
$$
Then almost surely $Y_0(0)<0<Y_0(1)$.

We use the Poisson invariance statement in \cite[Section 1, immediately after (1.1)]{FSW15}: for Brownian particles ordered increasingly and reflected upwards from their left neighbours, a homogeneous Poisson initial point process remains Poisson at each fixed time.
Their labeling $x_0(0)\le0<x_1(0)$ agrees with $Y_0(0)<0<Y_0(1)$ here.
Their intensity-one normalization gives intensity $\lambda$ by the scaling $x_k(t)\mapsto\lambda^{-1}x_k(\lambda^2t)$.

We identify the level boundaries of $h_t$ with this system.
Since $h_0(x)\to-\infty$ as $x\to-\infty$, the evolution $h_t$ cannot be directly linked to reflected Brownian motions.
For $m\ge1$, we introduce the truncation $h_0^{(m)}=h_0\vee(-m)$ and we let
$h_t^{(m)}$ be its evolution under \eqref{defevolution}.
By Lemma \ref{lemma:evolutionboundaries}, applied after adding
$m$ to the height, its level boundaries form the semi-infinite
reflected Brownian system obtained by retaining the initial
particles with labels $k>-m$, letting the leftmost particle
move freely and all further particles be reflected from their left neighbour.

By \cite[Proposition 2.4]{WFS17}, these truncated systems
converge, jointly for finitely many labels, to the infinite
reflected Brownian system. 
To verify admissibility, fix $k\in\Z$. For all sufficiently large $m$,
$Y_0(k)-Y_0(-m)$ is bounded below by
\begin{equation}\label{sumofexp}
Y_0(\min\{k,0\})-Y_0(-m)=\sum_{j=1}^{m+\min\{k,0\}}\xi_j
\end{equation}
where $\xi_j$ are independent exponential random variables of rate $\lambda$.
The lower bound by \eqref{sumofexp} and exponential Markov inequality therefore gives
$$
P\left(Y_0(k)-Y_0(-m)\le m^{3/4}\right)
\le P\left(e^{-\lambda\sum_{j=1}^{m+\min\{k,0\}}\xi_j}\ge e^{-\lambda m^{3/4}}\right)
\le e^{\lambda m^{3/4}}2^{-m-\min\{k,0\}}.
$$
The right-hand side is summable in $m$, so the admissibility
condition in \cite[Definition 2.2]{WFS17} holds with $\chi=3/4$.

On the other hand, directly from \eqref{defevolution},
$$
h_t^{(m)}(y)\downarrow h_t(y).
$$
The limit is defined by \eqref{defevolution} and it is finite by Theorem \ref{t:large-dev},
since the quadratic growth of $\mathcal B(x,0;y,t)$ as $x\to-\infty$
dominates the linear growth of $-h_0(x)$, while as $x\to\infty$, $h_0(x)\to\infty$ and the distance is nonnegative.

Therefore, for each $k\in\Z$,
$$
\sup\{y:h_t^{(m)}(y)<k\}
\uparrow
\sup\{y:h_t(y)<k\}.
$$
Thus the level boundaries of $h_t$ have the joint law of
the infinite reflected Brownian particle system.
At each fixed time their point process is therefore Poisson
with intensity $\lambda$. Normalizing the  counting
function at zero proves the claim; simplicity of the Poisson process ensures that its jumps have size one.
\end{proof}

\section{Dynamics}\label{s:dynamics}

We prove Theorem \ref{thm:fdconvdyn} in this section.
The proof proceeds similarly to that of Theorem 1.4 in \cite{VV26}.
We first introduce the left neighbourhoods of the rescaled random walk web distance as well as those of the lower semicontinuous Brownian web distance.
For any fixed non-negative integer $j$, time $s\in\R$, interval $I\subseteq\R$ and dynamical time $\tau\ge0$, we define
\begin{align}
R_{j,n}(I,s;\tau)&=\{(x,t):\inf_{z\in I}D^\RW_n(x,t;z,s;\tau)\le j\}\label{defRj}\\
Q_j(I,s;\tau)&=\{(x,t):\inf_{z\in I}D^\Br(x,t;z,s;\tau)\le j\}\label{defQj}
\end{align}
and we let $r_{j,n}^\pm$ and $\rho_j^\pm$ the boundary curves of the regions $R_{j,n}(I,s;\tau)$ and $Q_j(I,s;\tau)$ respectively.
We mention that the boundary curves can be written recursively using Skorokhod reflection of a random walk or a Brownian motion on the previous curve,
see formulas (3.2), (3.4) and Proposition 3.5 in \cite{VV26}.
The next analog of Proposition 4.2 in \cite{VV26} in the dynamical setting shows that boundary curves converge.

\begin{proposition}\label{prop:boundaries}
Let $0\le\sigma_1<\dots<\sigma_k$ be fixed.
Consider a countable collection of $s_i,I_i,\tau_i$ where $s_i\in\R$, $I_i\subseteq\R$ is a closed interval
and $\tau_i\in\{\sigma_1,\dots,\sigma_k\}$ for all $i$.
There exists a coupling under which, for all $i$, the boundary curves $r_{j,n}^\pm$ of the region $R_{j,n}(I_i,s_i;\tau_i)$ defined in \eqref{defRj}
converge in probability to the boundary curves $\rho_j^\pm$ of the region $Q_j(I_i,s_i;\tau_i)$ defined in \eqref{defQj} uniformly on compact sets.
\end{proposition}

The proof of Proposition \ref{prop:boundaries} relies on the following lemma.

\begin{lemma}\label{lemma:fdconvwebanddual}
Let $(Y^{(n)}(\tau))_{\tau \ge 0}$ be the diffusively rescaled DyDW, and $(\wh{Y}^{(n)}(\tau))_{\tau \ge 0}$ be its dual. Let $(B(\tau))_{\tau \ge 0}$ be the DyBW, and $(\wh{B}(\tau))_{\tau \ge 0}$ its dual. For any deterministic times $0 \le \tau_1 <\dots <\tau_k <\infty$, we have
\begin{multline}
(Y^{(n)}(\tau_1), \dots, Y^{(n)}(\tau_k)) \times (\wh{Y}^{(n)}(\tau_1), \dots, \wh{Y}^{(n)}(\tau_k))\\
\stackrel{\d}\Rightarrow(B(\tau_1), \dots, B(\tau_k)) \times (\wh{B}(\tau_1), \dots, \wh{B}(\tau_k))
\end{multline}
with respect to the product topology on $\Hc^k \times \widehat{\Hc}^k$.
\end{lemma}

\begin{proof}[Proof of Proposition \ref{prop:boundaries}]
By the Skorokhod representation theorem, the convergence in distribution of Lemma \ref{lemma:fdconvwebanddual} can be realized by a coupling
where the random walk webs $Y^{(n)}(\tau_i)$ and their duals $\wh Y^{(n)}(\tau_i)$ for $i=1,\dots,k$ converge
to the Brownian webs $B(\tau_i)$ and duals $\wh B(\tau_i)$ almost surely uniformly on compact sets.
Using Donsker's invariance principle and the continuity of Skorokhod reflection,
the argument in the proof of Proposition 4.2 in \cite{VV26} implies that under the coupling above
$r_{j,n}^\pm\to\rho_j^\pm$ in probability uniformly on compact sets.
As in \cite{VV26}, this holds for countably many $s_i,I_i,\tau_i$ and $j$ simultaneously.
\end{proof}

\begin{proof}[Proof of Theorem \ref{thm:fdconvdyn}]
By Proposition \ref{prop:boundaries},
the random walk webs $Y^{(n)}(\tau_i)$, their duals and countably many region boundaries converge jointly in probability.
Since the product of $(\Hc\times\wh\Hc)^k$ with countably many copies of $\Pi_0$ is still separable,
Skorokhod's representation theorem implies that the random walk webs and their duals can be coupled to the Brownian webs and their duals
in a way that the countably many boundaries all converge almost surely.
We let the countable collection consist of all $(s_i,I_i,\tau_i)$ triples where $s_i$ are rational, the endpoints of the intervals $I_i$ are rational
and $\tau_i\in\{\tau_1,\dots,\tau_k\}$.
We fix the coupling given by Skorokhod's representation for the above choice of the countable collection.

The convergence in \eqref{fdconvRWtoBr} in the product topology on $(\mathcal E_*)^k$ means
that all coordinates converge in the epigraph sense under the coupling described above.
Each coordinate marginally is a sequence of rescaled random walk webs coupled to a Brownian web
such that all region boundaries at rational times and with target intervals with rational endpoints converge almost surely.
Hence the proof of Theorem 1.4 in \cite{VV26} applies and yields that
\begin{equation}
\epi D^\RW_n(\tau_i)\to\epi D^\Br(\tau_i)
\end{equation}
almost surely in $\mathcal E_*$ for all $i=1,\dots,k$ which proves \eqref{fdconvRWtoBr}.
Finally, \eqref{fdconvSJtosBr} follows by the coordinatewise application of Theorem \ref{thm:convtransform}
for all $k$ values of the dynamical time under the same coupling.
\end{proof}

\begin{proof}[Proof of Lemma \ref{lemma:fdconvwebanddual}]
For each $i$, we have that
\begin{equation}
(Y^{(n)}(\tau_i), \wh{Y}^{(n)}(\tau_i)) \stackrel{\d}\Rightarrow (B(\tau_i), \wh{B}(\tau_i)),
\end{equation}
and so it follows that the sequence $[(Y^{(n)}(\tau_1), \dots, Y^{(n)}(\tau_k)) \times (\wh{Y}^{(n)}(\tau_1), \dots, \wh{Y}^{(n)}(\tau_k))]_{n \ge 1}$ is tight.
Assume that we have convergence in distribution along the subsequence $(n_j)_{j \ge 1}$ to $(W(\tau_1), \dots, W(\tau_k)) \times (\wh{W}(\tau_1), \dots, \wh{W}(\tau_k))$.
Using the Skorokhod representation theorem, we may assume this convergence holds almost surely on some common probability space. 
By Proposition 3.3 of \cite{RS24}, the law of $(W(\tau_1), \dots, W(\tau_k))$ is the appropriate finite-dimensional distribution of the dynamical Brownian web.

If $(\mu_n)_{n \ge 1}$ are probability measures on $\Hc \times \wh{\Hc}$ which are supported on the set of pairs $(C,\wh{C})$
where each path of $\wh{C}$ does not cross any path of $C$, and $\mu_n \stackrel{\d}\Rightarrow \mu$, then $\mu$ must also be supported on this set.
It follows that for each $i$, no path from $\wh{W}(\tau_i)$ crosses any path from $W(\tau_i)$.

Next, we have that almost surely, for each $i$ and each $(x,t) \in \Q^2$, there is a unique continuous backwards path that starts at $(x,t)$ and does not cross any path of $W(\tau_i)$.
Indeed, this path can be defined for $s<t$ by
\begin{equation}
\wh{\pi}_{(x,t)}^{\tau_i}(s)=\inf\{y\in \R: \text{all paths started at }(y,s) \text{ in } W(\tau_i) \text{ are} \ge x \text{ at time } t\}.
\end{equation}
For a continuous backwards path, this is the only way to define it so that it will not cross any forward path, which shows uniqueness.
See \cite[Section 2.2]{TW98} for more details, including a proof that this definition does indeed yield a continuous path.

Thus, each $\wh{W}(\tau_i)$ must contain all those paths (for the time $\tau_i$), so it contains the closure of the set of all such paths.
Since $\wh{Y}^{(n)}(\tau_i) \stackrel{\d}\Rightarrow \wh{B}(\tau_i)$, the law of $\wh{W}(\tau_i)$ is the same as that of $\wh{B}(\tau_i)$.
By Proposition 5.7 of \cite{FINR04}, $\wh{B}(\tau_i)$ is the closure of the set of backwards paths started at $(x,t)$, for $(x,t) \in \Q^2$.
It follows that $\wh{W}(\tau_i)$ cannot contain anything else.
Indeed, if random closed sets $K\subseteq L$ have the same law, then for each open set $U$ in a countable basis, the indicators that $K$ and $L$ meet $U$ are ordered and have the same expectation, hence agree almost surely.
Simultaneously for this basis they agree, which implies $K=L$.
We conclude that $(\wh{W}(\tau_1), \dots, \wh{W}(\tau_k))$ must be the duals of $(W(\tau_1), \dots, W(\tau_k))$.
So every subsequential limit is the law of $(B(\tau_1), \dots, B(\tau_k)) \times (\wh{B}(\tau_1), \dots, \wh{B}(\tau_k))$, and this concludes the proof.
\end{proof}

\section{Open questions and problems}

\begin{question}
What is the limit of the (conditional) uniform law on geodesics between two points in the Sepp\"al\"ainen--Johansson model?
Is this random distribution a function of the Brownian web?
\end{question}

\begin{conjecture}[Universality]
Let $\mu$ be a nondegenerate probability  distribution on $\mathbb R_+$ so that for some $v>0$,
$\mu$ restricted to $[0,v)$ is a single point mass at $0$.
Consider the directed first passage model on $\Z^2$ in which i.i.d.\ weights with distribution $\mu$ are assigned to the horizontal edges
and $0$ to the vertical edges.
The model has a critical direction above which the rescaled passage time is almost surely $0$ and it is positive below this direction.
The scaling limit of the passage time under $0:1:2$ scaling along the critical direction is the Brownian web landscape.
\end{conjecture}

\begin{conjecture}[Distances in oriented percolation]
Consider directed first passage percolation with Bernoulli$(p)$ weights on all edges.
Let $\Theta$ be the set of slopes $\theta$ that there is an infinite  open (all zero) path with slope $\theta$.
Then in the directions $\inf \Theta$ and $\sup \Theta$, the first passage values converge to the Brownian web landscape.
See \cite{SS13} for the connection between oriented percolation and the Brownian web, and \cite{SS19} for the dynamical version.
\end{conjecture}

\begin{problem}
Show that the dynamical Sepp\"al\"ainen--Johansson first passage percolation distance converges as a c\`adl\`ag process,
rather than only in the sense of finite dimensional distributions.
First prove the existence of the limit, i.e.\ that of the dynamical Brownian web distance as a c\`adl\`ag process.
\end{problem}

\begin{problem}
Find the asymptotics of the number of geodesics in the Sepp\"al\"ainen--Johansson model
as it converges to its Brownian web landscape limit.
\end{problem}

\begin{problem}
The hitting time formula in Theorem \ref{thm:RBM-formula} involves a Poisson process.
Does it have a probabilistic interpretation?
Is it related to the stationarity of the Poisson counting process?
\end{problem}

\paragraph{AI use.}
The mathematical content of the paper was developed by the authors.
Claude Opus 5.5 and OpenAI ChatGPT Astra 6.1  were used for literature search, proofreading, creating figures, and minor edits of the final version. 
The authors take full responsibility for the content of the paper.

\paragraph{Acknowledgments.}
We are grateful to the Erd\H os Center and the R\'enyi Institute for the Focused workshop on the Brownian web, the Brownian net and their geometry.
The work of J.~Ransford was supported by an NSERC Postdoctoral Fellowship.
The work of B.\ Vet\H o was supported by the NKFI (National Research, Development and Innovation Office)
grants FK142124 and KKP144059 ``Fractal geometry and applications''.
The work of B.~Vir\'ag was supported by the NSERC Discovery Grant.

\bibliographystyle{dcu}
\end{document}